\documentclass[10pt,twoside]{siamart1116}

\usepackage[english]{babel}
\usepackage{graphicx,epstopdf,epsfig}
\usepackage{amsfonts,epsfig,fancyhdr,graphics, hyperref,amsmath,amssymb}

\usepackage{relsize}
\usepackage{paralist}

\usepackage[normalem]{ulem}
\newcommand{\msout}[1]{\ifmmode\text{\sout{\ensuremath{#1}}}\else\sout{#1}\fi}

\newcommand{\cancel}[1]{
	\ifmmode
	{\color{red}\msout{#1}}
	\else
	{\color{red}\sout{#1}}
	\fi
}

\newcommand{\replace}[2]{
	\ifmmode
	{\color{red}\msout{#1}}{\color{blue}\uwave{#2}}
	\else
	{\color{red}\sout{#1}}{\color{blue}\uwave{#2}}
	\fi
}

\newcommand{\HE}{}
\newcommand{\DoS}{}
\newcommand{\DoA}{}
\newcommand{\CA}{Seonhwa Kim}

\newcommand{\Names}{Yunhi Cho and Seonhwa Kim}
\newcommand{\Title}{Volume of Hypercubes Clipped by Hyperplanes\\ and Combinatorial Identities}

\newtheorem{remark}[theorem]{Remark}
\newtheorem{example}[theorem]{Example}

\newcommand{\reals}{\mathbb{R}}

\renewcommand{\theequation}{\thesection.\arabic{equation}}
\renewtheorem{theorem}{Theorem}[section]

\newtheorem*{LawAs*}{Lawrence's two conditions}

\def\Acal{\mathcal{A}}
\def\Hcal{\mathcal{H}}
\def\Ical{\mathcal{I}}

\def\abf{\mathbf{a}}
\def\bbf{\mathbf{b}}
\def\cbf{\mathbf{c}}
\def\ebf{\mathbf{e}}
\def\obf{\mathbf{o}}
\def\qbf{\mathbf{q}}
\def\vbf{\mathbf{v}}
\def\wbf{\mathbf{w}}
\def\xbf{\mathbf{x}}

\DeclareMathOperator{\vol}{vol}

\makeatletter
\providecommand*{\cupdot}{%
	\mathbin{%
		\mathpalette\@cupdot{}%
	}%
}
\newcommand*{\@cupdot}[2]{%
	\ooalign{%
		$\m@th#1\cup$\cr
		\sbox0{$#1\cup$}%
		\dimen@=\ht0 %
		\sbox0{$\m@th#1\cdot$}%
		\advance\dimen@ by -\ht0 %
		\dimen@=.5\dimen@
		\hidewidth\raise\dimen@\box0\hidewidth
	}%
}

\providecommand*{\bigcupdot}{%
	\mathop{%
		\vphantom{\bigcup}%
		\mathpalette\@bigcupdot{}%
	}%
}
\newcommand*{\@bigcupdot}[2]{%
	\ooalign{%
		$\m@th#1\bigcup$\cr
		\sbox0{$#1\bigcup$}%
		\dimen@=\ht0 %
		\advance\dimen@ by -\dp0 %
		\sbox0{\scalebox{2}{$\m@th#1\cdot$}}%
		\advance\dimen@ by -\ht0 %
		\dimen@=.5\dimen@
		\hidewidth\raise\dimen@\box0\hidewidth
	}%
}
\makeatother

\newcommand{\zerov}{{0_{\vbf}}}
\newcommand{\onev}{{1_{\vbf}}}
\newcommand{\bulletv}{{\bullet_{\vbf}}}
\newcommand{\starv}{{*_{\vbf}}}

\begin{document}

\bibliographystyle{plain}

\setcounter{page}{1}

\thispagestyle{empty}

 \title{\Title\thanks{Received
 by the editors on \DoS.
 Accepted for publication on \DoA. 
 Handling Editor: \HE. Corresponding Author: \CA}}

\author{
%
%
%
%
Yunhi Cho\thanks{Department of Mathematics, University of Seoul, Seoul 02504 , Korea
	(yhcho@uos.ac.kr). Supported by the 2019 Research Fund of University of Seoul.}
%
\and
Seonhwa Kim\thanks{
	School of Mathematics, Korea Institute for Advanced Study, Seoul 02455, Korea
	(seonhwa17kim@kias.re.kr, @gmail.com). Supported by 
	 the National Research Foundation of Korea (NRF) grant funded by the Korea government (MSIT) (No. 2019R1C1C1003383), and supported by the Institute for Basic Science (IBS-R003-D1).}
}

\markboth{\Names}{\Title}

\maketitle

\begin{abstract}
	There is an elegant expression for the volume of hypercube $[0,1]^n$ clipped by a single hyperplane.  In the article the formula is generalized to the case of more than one hyperplane. An important foundation for the result is Lawrence's formula and  a way to weaken two restrictions of simplicity and non-parallelness in his formula is also considered.  Several concrete  volume formulas of clipped hypercubes  are derived explicitly  and the corresponding  combinatorial identities are obtained as an application.
\end{abstract}




\section{Introduction}
A  unit hypercube is a convex polytope defined by $[0,1]^n$ in $\reals^n$. It may be  a very  basic geometric object and the simplest convex polytope, but it still has  interesting unsolved questions 
(for examples, see \cite{zong_what_2005}). 
It has turned out that the computation of  the volume of convex polytopes is  algorithmically hard \cite{dyer_complexity_1988} since it usually requires difficult work like vertex/facet enumeration, even for the case of hypercubes clipped by only one hyperplane \cite{khachiyan_problem_1989}.
For these reasons,  both numerical approximations and exact calculations  have been extensively studied from an algorithmic point of view.
%

In this paper, 
we will only focus  on  closed and  concrete formulas  using matrix computations. 
Although this doesn't require  heavy machinery (but is technically complicated),
the resulting volume formula is fairly concrete and hence has some applications like producing  a certain  class of combinatorial identities.  Furthermore, 
it would be interesting to investigate the relation between our formulation focused on
$[0,1]^n$ and Filliman's study \cite[Section 3]{filliman_volume_1992}.


For the easiest case of a hypercube clipped by only one hyperplane, there is an interesting  simple  formula giving the volume as the following. The notation $ |\zerov|$  for a vector $\vbf \in \reals^n$ indicates the number of zeros in the entries and $F^0$ denotes the set of vertices of $[0,1]^n$. The half space  $H_1^+$ is given by 
$$\{\xbf ~|~ g_1(\xbf):= \abf \cdot \xbf +r_1 = a_1x_1+ a_2x_2+ \cdots + a_n x_n  + r_1  \geq0\}$$
with  $\prod_{t=1}^n  a_t \neq 0$. Then we have
\begin{theorem}\label{thm:monthly}	
	\begin{equation*} 
	\vol ( [0,1]^n \cap H_1^+ )	= \sum_{ \vbf \in F^0 \cap  H_1^+}
	\frac{(-1)^{|\zerov|} g_1(\vbf)^n} {n!\prod_{t=1}^n  a_t}.
	\end{equation*} 
\end{theorem}
This formula seems to have first appeared in  \cite{barrow_spline_1979}, but a very similar idea seems to go back
much earlier \cite{polya_berechnung_1913}.  Although it  has been revisited several times (for examples, see Section 2 in \cite{marichal_slices_2008}), as far as the authors know,  a volume formula for the case of more than one hyperplane had not been seriously studied  yet. 
We generalize this formula to the case of an arbitrary number of hyperplanes on the basis of Lawrence's work \cite{lawrence_polytope_1991}.

Our volume formulas have concrete and explicit expressions, which have  some benefits  for the case of a small number of hyperplanes but a large dimensional cube. One can use our formula even for the case of sufficiently many  hyperplanes making up a fully general polytope. But the greater the number of hyperplanes, the less useful our formula seems to be, because the characteristics coming from the shape of cube tend to disappear and the formula essentially becomes almost the same as Lawrence's.

The general formulas will be described in Section \ref{sec:severalversion}.
Before doing that, let's look at two hyperplanes, which is a corollary of Theorem \ref{thm:mainvee} (the detailed description is given in Section \ref{sec:twoplanes}).
Let the  half spaces  $H_1^+$ and $H_2^+$ be  given by 
\begin{align*}
H^+_1 &=\{\xbf ~|~ g_1(\xbf) := a_1x_1+a_2x_2+\cdots+a_n x_n+r_1 \geq0\} \text{ and }\\
H^+_2 &=\{ \xbf ~|~ g_2(\xbf):= b_1x_1+b_2x_2+\cdots+b_nx_n+r_2\geq0\}
\end{align*} 
with good clipping conditions. 	
See Section \ref{sec:notation} for the notation and Section \ref{sec:goodclipping} for the definition of good clipping conditions. Then  Theorem \ref{thm:monthly} is generalized to the case of two hyperplanes as follows.
\begin{corollary}\label{thm:twoplanes}
	\begin{align*}
	\vol([0,1]^n \cap H_1^{+} \cap H_2^{+})
	=&\sum_{\vbf\in F^0 \cap H_1^{+} \cap H_2^{+}}
	\frac{(-1)^{|\zerov|}g_2(\vbf)^n}{n! \prod_{t=1}^n b_t}\\
	-&\sum_{\vbf\in F^1 \cap H_1 \cap H_2^{+}}
	\frac{(-1)^{|\zerov|}~a_{{*(\vbf)}}^{n}~g_2(\vbf)^n}
	{\;\; \;  n! ~ |a_{*(\vbf)}|~b_{*(\vbf)} 	\prod\limits_{t\in[n]\setminus *(\vbf) } \left| \begin{array}{cc} a_{{*(\vbf)}} & b_{*(\vbf)} \\ a_{t} & b_{t} \\ \end{array} \right|}.
	\end{align*}

\end{corollary}

Interestingly, 
volumes of hypercubes  clipped by various choices of hyperplanes  produce  non-trivial combinatorial identities. 
Let us see  several examples. 
\begin{theorem}\label{thm:generalizedRuiz}  For arbitrary $y \in \reals $, $a_1,a_2,\dots,a_n \in \reals$  and an integer $n \geq 0$, 
	$$y^n+\sum_{i=1}^n \;\; \sum_{1\leq t_1 < t_2 <\cdots < t_i \leq n} (-1)^{i}(y+a_{t_1}+\cdots+a_{t_i})^n  = (-1)^n n! a_1 a_2 \cdots a_n$$
	or equivalently,
	$$\sum_{i=1}^n \;\; \sum_{1\leq t_1 < t_2 <\cdots < t_i \leq n} (-1)^{i}(a_{t_1}+\cdots+a_{t_i})^k  = 
	\begin{cases}
	-1, & \textup{if }\;\; k=0 \\
	0, & \textup{if }\;\; k=1,\dots,n-1 \\
	(-1)^n n! a_1 a_2 \cdots a_n, & \textup{if }\;\; k=n.
	\end{cases}
	$$
\end{theorem}
This is related to an old Prouhet-Tarry-Escott problem \cite{lehmer_tarry-escott_1947} and there are many  not so difficult proofs. But it can be proved by  Theorem \ref{thm:monthly} and we think the geometric proof using a clipped cube  is a new approach.

The next formula also may be obvious to someone familiar with  Vandermonde matrices or Lagrange's interpolation formula.
However, it can also be proved directly by  Corollary~\ref{thm:twoplanes} whose geometric description is  a simplex clipped by one hyperplane (see \cite{cho_volume_2001}) or the hypercube $[0,1]^n$ clipped by two hyperplanes: a general hyperplane and a special hyperplane which passes through all standard basis vectors.

\begin{theorem}\label{thm:clippedsimplex} 
	For arbitrary $y \in \reals $, distinct non-zero $a_1,a_2,\dots,a_n \in \reals$ and an integer $n \geq 0$, 	
	\begin{equation*}\label{e3}		
	\frac{y^n}{a_1a_2\cdots a_n} - \sum_{i=1}^n \frac {(y+a_i)^n}{a_i(a_1-a_i)(a_2-a_i)\cdots\widehat{(a_i-a_i)} \cdots(a_n-a_i)} = (-1)^n \end{equation*}
	or equivalently
	$$\sum_{i=1}^n \frac {a_i^k}{(a_1-a_i)(a_2-a_i)\cdots \widehat{(a_i - a_i)} \cdots (a_n-a_i)}=
	\begin{cases}
	\frac{1 }{a_1a_2\cdots a_n}, & \textup{if }\;\; k=-1 \\
	0,   & \textup{if }\;\; k=0,1,\ldots,n-2 \\
	(-1)^{n-1},   & \textup{if }\;\; k=n-1.
	\end{cases}
	$$
	where $\widehat{~~}$ means  omitting the term.
\end{theorem}
Furthermore, if we take $a_1=a_2=\dots=a_n=1$ in Theorem \ref{thm:generalizedRuiz}, then we get the following corollary which is a studied form in combinatorial enumeration (for example, \cite{ruiz_an_1996}).
\begin{corollary} For arbitrary $y \in \reals$ and an integer $n \geq 0$,
	$$\sum_{i=0}^n(-1)^{i}\binom{n}{i}(y+i)^n=
	(-1)^n n!
	$$
	or equivalently
	$$\sum_{i=0}^n(-1)^{i}\binom{n}{i}i^k=
	\begin{cases}
	0,   & \textup{if }\;\; k=0,1,\ldots,n-1 \\
	(-1)^n n!, & \textup{if }\;\; k=n.
	\end{cases}
	$$
\end{corollary}

Interestingly, all the above identities are unified under one umbrella via a volume expression for a particular clipped hypercube.
Before showing this, we introduce the following set-theoretic notation for the sake of convenience,
\begin{align*}
A:=\{a_1, a_2,\ldots, a_n\},~~~  \| A \|:=\sum\limits_{a\in A}a,  ~~~ A! :=\prod\limits_{a\in A}a, \\
R_A(a):=\prod\limits_{b\in A\setminus a}\frac {b}{b-a},~~~~~~~~~~~~~   
R_A(I) := \sum\limits_{a\in I} R_A(a).
\end{align*}
Then Theorem \ref{thm:generalizedRuiz} and Theorem \ref{thm:clippedsimplex} are written in an economic way as follows.
\begin{align*}
\sum_{I\subset A} (-1)^{|I|}\|I\|^k  &=
\begin{cases}
0, & \textup{if }\;\; k=0,1,\dots,n-1 \\
(-1)^n n!A!, & \textup{if }\;\; k=n,
\end{cases}
\end{align*}

\begin{align*}
\sum_{a\in A} R_A(a)a^k &= 
\begin{cases}
1,  & \textup{if }\;\; k=0 \\
0,   & \textup{if }\;\; k=1,\ldots,n-1 \\
(-1)^{n-1} A!,    & \textup{if }\;\; k=n.
\end{cases}
\end{align*}

We can also obtain the following identity which can be derived from Corollary \ref{thm:twoplanes} (for the proof, see Theorem \ref{thm:finalid}). 

\begin{theorem}\label{thm:finalcom}
	For $A=\{a_1,\dots,a_n\}$ and  an integer $l = 1,\dots,n$,
	\begin{align*}
	&\sum_{\substack{I\subset A\\|I|<{l}}}  {(-1)^{|I|} \|I\|^k}
	+  {\sum_{\substack{I\subset A\\|I|={l}}}}  (-1)^{l} \|I\|^k
	R_A(I)\\ 
	&~~~~~~~~~~~~~~~~=\begin{cases}
	0,   & \textup{if }\;\; k=0,1,\ldots,n-1 \\
	A! \sum\limits_{i=0}^{{l}-1} (-1)^{n-i} \binom{n}{i}({l}-i)^n   , & \textup{if }\;\; k=n.
	\end{cases} 
	\end{align*}
\end{theorem}	
The identity in  Theorem \ref{thm:finalcom} itself  may already be known, but the proof using the volume of  clipped hypercubes seems new. 
Finally we would like to remark that the above identities are all  symmetric functions. In Section \ref{sec:combidentity} and the Appendix, we will give several identities, some  symmetric and  others   not.

Let us outline our article. We will introduce  notation in Section \ref{sec:notation}, and review and reorganize Lawrence's formula and explain our $\epsilon$-perturbation method in Section \ref{sec:reviewLawrence}. The statements of main theorems and proofs will be given in Section \ref{sec:mainproof}. Several concrete examples will be presented with more explicit expressions  in Section \ref{sec:concereteexample}. 
In the final Section \ref{sec:combidentity}, we will derive a family of  combinatorial identities using the volume of clipped hypercubes.

\section{Notation }\label{sec:notation}

In this paper, the letters  $n$ and $m$ correspond to  the dimension of $\reals^n$ and  the number of hyperplanes respectively unless otherwise specified.  A single bold letter always denotes a vector in $\reals^n$ like $ \xbf = (x_1,\dots,x_n)$ and we abuse notation for column vectors and row vectors if it is not confusing. Let $\ebf_i$ denote the $i$-th vector in the standard basis of $\reals^n$. 

Let $K$ be the natural cell structure of   unit hypercube   $[0,1]^n$ in  $\reals^n$ and let $K^d$ denote its $d$-skeleton. We define the \emph{open $d$-skeleton} $F^d$ as $K^d \setminus K^{d-1}$. Then, 
\begin{align}\label{eqn:openskel}
[0,1]^n = \bigcup_{d=0}^n K^d = \bigcupdot_{d=0}^n F^d ,
\end{align}
where the $\cupdot$ symbol denotes disjoint union. For example,
$[0,1]^2$ consists of
{four} points $F^0$,
{four} open intervals $F^1$ and
{one} open rectangle $F^2$.

\subsection{Index manipulation}
Let $[n]$ denote the \emph{ordered set} $\{1,2,\ldots,n\}$ that is an index set for the standard basis of $\reals^n$.
We will use ordered sets for indices because the sign of a minor of a matrix  is sensitive to the order of indices. 
Let $A_I^J$ and $(A)_I^J$ denote  a minor and a submatrix with indices $I$ and $J$ of a matrix $A=(a_{i,j})$ respectively.

For example, let $I = \{1,3\}$ and $J=\{2,4\}$. Then 
$$(A)_I^J=\left[ \begin{array}{cc} a_{1,2} & a_{1,4} \\ a_{3,2} & a_{3,4} \\ \end{array} \right] \text{ and } A_I^J = \left| \begin{array}{cc} a_{1,2} & a_{1,4} \\ a_{3,2} & a_{3,4} \\ \end{array} \right| = \det
\left[ \begin{array}{cc} a_{1,2} & a_{1,4} \\ a_{3,2} & a_{3,4} \\ \end{array} \right] . $$ 


Let an ordered set $I=\{i_1,i_2,\ldots,i_s\} \subset [n]$. Elementary arithmetic operations  with an ordered set and a number are done  entrywise, for example $2I-1 = \{2 i_1 -1 ,\ldots,2 i_s -1 \}$. 
We call an index $I$  \emph{well-ordered} if $i_1 < i_2< \cdots <i_s$. We consider two different notions of union operation for ordered sets.  One is the \emph{ordered union} $\cup$ respecting the order between two well-ordered indices, for instance, for $t \not \in I$,
$$I\cup \{t\} \;\;:=\;\; \{i_1,i_2,\ldots,t,\ldots,i_s\} \;\;\textrm{ when } i_1<i_2<\cdots<t<\cdots<i_s .$$
The other is  the \emph{joining union} $\vee$ as concatenation as follows,
$$I\vee \{t\} \;\;:=\;\; \{i_1,i_2,\ldots,i_s,t\}.$$ 
We remark that the joining union is defined no matter whether the constituent sets are well-ordered or not, but the ordered union is defined only for well-ordered sets. In general, the result of a joining union is not well-ordered and might be an ordered multi-set.

We abbreviate  a set of one element $\{x\}$ to $x$  omitting the brace symbols, for example, $I \vee \{t\} =: I \vee t$.
Let $I$ and $J$ be two ordered sets consisting of the same elements. 
Then $\sigma({I,J})$ denotes  the parity of the permutation between the two ordered sets $I$ and $J$, for example $\sigma(a\vee b, b \vee a )=-1$.

Let $| \cdot |$ and $|| \cdot ||$ denote the cardinality  and the total sum of elements of a given set respectively. We remark that  $\|\varnothing \|^k = 0^k= 1$ when $k=0$.
For $\vbf=(v_1,v_2,\dots,v_n)$ in $\reals^n$,
we define  the notation $\zerov$, $\onev$, $\bulletv$ and $\starv$ 
which denote 
ordered sets of indices satisfying the following. 
\begin{align}\label{def:v01}
\begin{aligned}
\zerov &:=\{i \in [n] \mid v_i=0 \}, 
~~&\onev &:=\{i \in [n] \mid v_i=1 \}, \\
\starv &:=\{i \in [n] \mid v_i \neq 0,1 \}, 	 	 
~~&\bulletv &:= \zerov \cup \onev = [n]     \setminus  \starv .
\end{aligned}
\end{align} 
In particular, we define  functions $*_i:\reals^n \rightarrow [n]$ and $\bullet_i : \reals^n \rightarrow [n]$ by indicating the $i$-th entry of  
$\starv$ and $\bulletv$ of increasing order respectively, i.e.,
\begin{align*}
\starv = \{i\in[n] \mid v_i\ne 0,1\} &=\{*_1(\vbf), *_2(\vbf),\ldots, *_{|\starv|}(\vbf) \},    \\
\bulletv = \{i\in[n] \mid v_i = 0,1\} & =\{\bullet_1(\vbf), \bullet_2(\vbf),\ldots, \bullet_{|\bulletv|}(\vbf)\}     .
\end{align*}

When we consider a set of only one element   then we omit the index letter like $ *(\vbf) := *_1(\vbf)$. 
To help understanding, let us see an example. Let  $\vbf=(0,1,\frac 13,0,0,\frac 35,1,\frac 18)\in F^3$, then we get 
\begin{align*}
&\starv=\{3,6,8\},\;\; \zerov=\{1,4,5\},\;\; \onev=\{2,7\}, \\
&|\starv|=3,\;\; |\zerov|=3,\;\;|\onev|=2, \\
&\|\starv\|=3+6+8=17,\;\;\|\zerov\|=1+4+5=10,\;\;\|\onev\|=2+7=9,\\
& *_1(\vbf)=3,\;\; *_2(\vbf)=6,\;\; *_3(\vbf)=8, \\
& \bullet_1(\vbf)=1,\;\; \bullet_2(\vbf)=2,\;\; \bullet_3(\vbf)=4,\;\; \bullet_3(\vbf)=5,\;\; \bullet_3(\vbf)=7. 
\end{align*}
Finally,  we remark that the following always holds by definition:
\begin{align*}
\zerov  \cup \onev \cup \starv &= \bulletv \cup \starv = [n]. 
\end{align*}

\subsection{Hyperplane matrices}\label{sec:hyperplaneNotation}
Throughout the article, hyperplanes and half spaces are given by  
$$H_i := \{ \xbf \;|\; g_i(\xbf)=0\} \textrm{ and } H^+_i:= \{ \xbf \;|\; g_i(\xbf)\geq0\},$$ where the linear coefficients are the following:
\begin{align*}
g_1(\xbf) &:= \abf_1\cdot \xbf + r_1 &=\;&  a_{11}x_1+a_{21}x_2+\cdots+a_{n1}x_n+r_1, \\ 
&&\vdots\\
g_{m-1}(\xbf) &:= \abf_{m-1} \cdot \xbf + r_{m-1} &=\;&  a_{1,m-1}x_1+a_{2,m-1}x_2+\cdots+a_{n,m-1}x_n+r_{m-1},\\ 
g_m(\xbf) &:= \abf_m \cdot \xbf + r_m &=\;&  a_{1,m} x_1+a_{2,m} x_2+\cdots+a_{n,m} x_n+r_m. 
\end{align*}
These coefficients form an $n \times m$ matrix
$ {{A}}$ as the following.
$$ {{A}}:=(\abf_1, \abf_2, \ldots,\abf_{m-1},\abf_{m})=\left[
\begin{array}{ccccc}
a_{1,1} & a_{1,2} &      &a_{1,k-1} &a_{1,m}\\
a_{2,1} & a_{2,2} &      &a_{2,k-1} &a_{2,m}\\
\vdots & \vdots &\cdots&\vdots    &\vdots\\
a_{n,1} & a_{n,2} &      &a_{n,k-1} &a_{n,m}\\
\end{array}
\right].$$
In particular we will take  the last $ g_m(\xbf)$ and  $H_m$ as the \emph{auxiliary function} and the \emph{auxiliary hyperplane} respectively. Let $H^+$  denote the intersection of all  the half spaces $H_i^+$,

$$H^+ = \bigcap_{i\in [m]} H_i^+.$$

Let $I$ be a set of indices for several hyperplanes usually not including the auxiliary plane, i.e. $I \subset [m-1]$ and let $H_I$ denote the intersection of   $H^+ \setminus H_m$ and the hyperplanes $H_i$ for $i\in I$, i.e.,
$$H_I:= \bigcap\limits_{i\in I} H_i \cap H^+ \setminus H_m.$$ 
We remark that  we  remove  the auxiliary plane $H_m$ from the definition of $H_I$ because we are going to ignore so-called  \emph {degenerate} vertices (see Section \ref{sec:degenerate}).

\section{A volume formula for convex polytopes}\label{sec:reviewLawrence}

\subsection{A review on volume computations}\label{sec:pictorialreview}
\begin{figure}[h!]\label{fig:decomp}
	\centering
	\def\svgwidth{120mm}
	\input{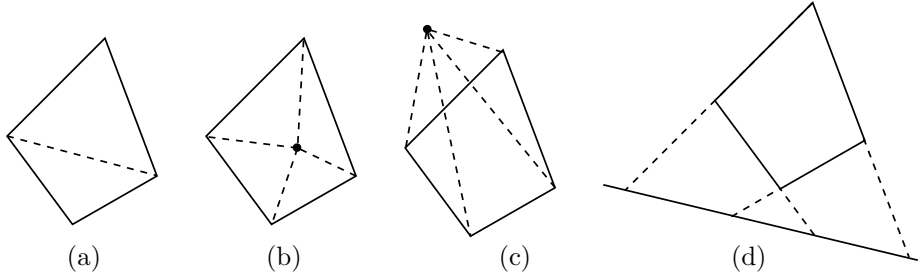}
	\caption{Typical decomposition methods for a convex polytope}
\end{figure}

Let us briefly review  conceptual methods to compute the exact volume of convex polytopes in a pictorial way. 
It is well known that the volume of an $n$-parallelotope and an $n$-simplex given by $\vbf_1,\vbf_2,\dots, \vbf_n$ in $\reals^n$ are  $|\det( \vbf_1 \vbf_2\dots \vbf_n)|$ and $\frac{1}{n!}|\det( \vbf_1\vbf_2\dots \vbf_n)|$ respectively.
An elementary strategy for computing the volume of a  polytope is to decompose the polytope into a signed summation of simplices. In fact,  many volume computing algorithms rely entirely on the method of decomposition as in Figure \ref{fig:decomp}. 
Case (a) is the most obvious decomposition that always exists  because of convexity.  Cases (b) and (c) are essentially the same except for the position of an auxiliary point. These decompositions are quite elementary but interesting because they make an identity between the volume of a polytope and the volumes of its facets (for examples, see J.B. Lasserre \cite{lasserre_analytical_1983}).
Case (d)  is, in some sense, a dual approach to (c) since it uses an auxiliary plane instead of an auxiliary point.  

A volume formula using decomposition (d) was first proposed by J. Lawrence \cite{lawrence_polytope_1991} and  is a base camp for this study. 
Let us review the results.
A convex compact polytope $P$ in Euclidean $n$-space is given by 
$$P=\bigcap_{i\in [m]} H_i^+\subset \reals^n.$$ 
Consider an auxiliary function $f(\xbf)$ with a hyperplane $f(\xbf)=0$.  Then we need two requirements as follows.

\begin{LawAs*} \leavevmode \vspace{-0.3em}
	\begin{description}
		\item[simplicity] All vertices of a polytope $P$ are  \emph{simple}, which means that the degrees of vertices in $P$ are the same as the dimension of  $P$.
		\item[non-parallelness]  The auxiliary function $f$ is non-constant on each edge of $P$.	
	\end{description}
	
\end{LawAs*} 
Now, we have the following volume formula.

\begin{theorem*}[J. Lawrence] 
	If  $P$ satisfies  \emph{ Lawrence's  condition}, the volume  is computed by
	\begin{equation}\label{eq:Lawrence}
	\vol (P)=\sum_{\vbf: \text{ a vertex of }P} N_\vbf
	\quad\textup{ with }\quad N_\vbf=\frac{f(\vbf)^n}{n!\delta_\vbf\gamma_1 \gamma_2 \cdots \gamma_n}
	\end{equation}
	where $\delta_\vbf$ is the absolute value of the determinant of the $n\times n$ matrix whose columns are $-\abf_{i_1}, -\abf_{i_2},\ldots,-\abf_{i_n}$. Here, a vertex $\vbf$ of $P$ is the intersection of $H_{i_1},\dots,H_{i_n}$ with $i_1<i_2<\dots<i_n$ and $\gamma_1, \gamma_2, \dots, \gamma_n$ are obtained by $\abf_m=-(\gamma_1 \abf_{i_1}+\cdots +\gamma_n\abf_{i_n})$. 
\end{theorem*}

For a more detailed description, see J. Lawrence \cite{lawrence_polytope_1991}.
From a geometric point of view, each $N_\vbf$ corresponds to a signed volume of a \emph{simplex}  which is \emph{projected} from a vertex $\vbf$ to the plane $\{ \xbf \in\reals^n |~f(\xbf)=0\}$ as in Case (d) in figure \ref{fig:decomp}. Note that the projected simplex possibly has infinite volume, but non-parallelness guarantees that all vertices of the simplex are located in a bounded region of $\reals^n$, not at infinity)

Note that the non-parallel condition is a sufficient condition but not a necessary condition for the volume of a projected simplex not to be $\infty$. By observation we found a case that violates the requirements, where we can  exploit Lawrence's formula without an extra effort.  Before doing that, we review his formula  concretely in the following section.

\subsection{An explicit form of Lawrence's formula}
Lawrence's formula  has sometimes been expressed in a quite explicit form in the literature (for example, see  \cite[p.393]{gritzmann_complexity_1994}). But the known formulations are still not enough  to proceed for our purpose. 
So we rewrite his formula in a closed form explicitly with linear coefficients.
From now on, we put an auxiliary plane as the last hyperplane $H_m$. In this section, $H_m$ is redundant, i.e. $P=P \cap H_m^{+}$ and $H_m\neq H_i$ for any $i\in[m-1]$.
We remark that $m-1 > n$ holds in order to make a convex compact polytope. 

\begin{theorem}\label{thm:mainLawrence}
	Let a convex polyhedron $P =  \bigcap\nolimits_{i\in[m-1]} H_i^+$ satisfy Lawrence's conditions.
	Then the volume is 
	\begin{align*}
	\vol (P )
	=\mathlarger\sum_{\substack{I\subset[m-1]\\|I|=n}}
	~~\mathlarger\sum_{\vbf\in H_I}
	\frac{(-1)^{ \frac{n(n+1)}{2}} (g_m(\vbf) {{A}}_{[n]}^{I})^n}{n!| {{A}}_{[n]}^{I}| \prod\limits_{t\in I}  {{A}}_{[n]}^{I{\setminus t \cup m}}}.
	\end{align*}
	
\end{theorem}
\begin{remark}
	In the formula, the second summation consists of either an empty summand or only one summand. In spite of the redundancy of the expression,  we persist in this inefficient form to be compatible with the main theorem on a clipped hypercube. We remark that
	$I{\setminus t \cup m} = I{\setminus t \vee m}$ because $m$ is the last element.
\end{remark}
\begin{proof}
	Let us just do a direct computation in our setup from Lawrence's formula  (\ref{eq:Lawrence}).
	Each vertex $\vbf$ is an  intersection of exactly $n$  hyperplanes other than the auxiliary plane, $$H_{i_1},H_{i_2},\ldots,H_{i_n}~ (1\leq i_1 < i_2 < \cdots < i_n \leq m-1)$$ with $\vbf=\bigcap_{t=1}^n H_{i_t}$.
	Let $I:=I_{\vbf}:=(i_1,\dots,i_n)$. Then 
	$$(- {{A}})_{[n]}^I=(-\abf_{i_1},-\abf_{i_2},\ldots,-\abf_{i_n}).$$
	
	Let  $\pmb{\gamma} := \pmb{\gamma}_\vbf := (\gamma_1,\ldots,\gamma_n)^t$.
	It follows from the definition in \cite{lawrence_polytope_1991} that  $\pmb{\gamma}$ is defined as satisfying  
	$$\abf_m=(- {{A}})_{[n]}^I \pmb{\gamma}.$$
	
	So we have $\pmb{\gamma}=((- {{A}})_{[n]}^I)^{-1}\abf_m$.
	
	Let $(x_{\mu \nu}) :=((- {{A}})_{[n]}^I)^{-1}$. Then  by Cramer's rule we get $x_{\mu \nu}=\frac{(-1)^{\mu +\nu}((- {{A}})^I_{[n]})_{\nu,\mu}}{\det((- {{A}})^I)}$ where $((-A)^I_{[n]})_{\nu,\mu}$ means  the $(\nu,\mu)$-minor of the submatrix $(-A)_{[n]}^I$, i.e. $|(-A)_{[n]\setminus \nu}^{I \setminus i_{\mu}}|$ with $I=\{i_1,i_2,\dots,i_n\}$.
	Then 
	
	$$\begin{aligned}\mathbf{\gamma}=\left[\begin{array}{c}
	\sum_j a_{j,m} x_{1,j}\\
	\sum_j a_{j,m} x_{2,j}\\
	\vdots\\
	\sum_j a_{j,m} x_{n,j}\\
	\end{array}\right]&=
	\frac{1}{\det((- {{A}})^I)}\left[\begin{array}{c}
	\sum_j (-1)^{1+j}a_{j,m} ((- {{A}})^I_{[n]})_{j,1}\\
	\sum_j (-1)^{2+j}a_{j,m} ((- {{A}})^I_{[n]})_{j,2}\\
	\vdots\\
	\sum_j (-1)^{n+j}a_{j,m} ((- {{A}})^I_{[n]})_{j,n}\\
	\end{array}\right]\\
	&=
	\frac{1}{\det((- {{A}})^I)}\left[\begin{array}{c}
	\det(\abf_m,-\abf_{i_2},\ldots,-\abf_{i_n})\\
	\det(-\abf_{i_1},\abf_m,\ldots,-\abf_{i_n})\\
	\vdots\\
	\det(-\abf_{i_1},-\abf_{i_2},\ldots,\abf_m)\\
	\end{array}\right].\end{aligned}$$ Hence
	$$\begin{aligned}
	\gamma_1\gamma_2\cdots\gamma_n&=
	\prod_{j=1}^{n}\frac{\det(-\abf_{i_1},\ldots,-\abf_{i_{j-1}},\abf_m,-\abf_{i_{j+1}},\ldots,-\abf_{i_n})}{\det((- {{A}})_{[n]}^I)}\\
	&=\prod_{j=1}^{n}\frac{(-1)^{n-1}\det(\abf_{i_1},\ldots,\abf_{i_{j-1}},\abf_m,\abf_{i_{j+1}},\ldots,\abf_{i_n})}{(-1)^n\det( {{A}}_{[n]}^I)}\\
	&=(-1)^n(-1)^{\frac{n(n-1)}{2}}\prod_{j=1}^{n}\frac{\det(\abf_{i_1},\ldots,\abf_{i_{j-1}},\abf_{i_{j+1}},\ldots,\abf_{i_n},\abf_m)}{\det( {{A}}_{[n]}^I)}\\
	&=\frac{(-1)^{\frac{n(n+1)}{2}}}{\det( {{A}}_{[n]}^I)^n}\prod_{t\in I}~ {{A}}_{[n]}^{I{\setminus t \cup m}}.
	\end{aligned}$$
Also we get
\begin{align*}
\delta_\vbf=|\det((- {{A}})_{[n]}^I)|=|\det( {{A}}_{[n]}^I)|=| {{A}}_{[n]}^{I}|.
\end{align*}
\end{proof}

\subsection{Degenerate vertices and edges}\label{sec:degenerate}
Let us consider the situation  of Theorem \ref{thm:mainLawrence} such that the last hyperplane $H_m$ is not redundant. In other words we choose an auxiliary plane $\{f(\xbf)=0\}$ as the hyperplane $H_m=\{\xbf \mid g_m(\xbf)=0\}$, when a polytope $P$ is given by a non-redundant intersection of $H_i$'s for $i \in [m]$, i.e., $P \neq  \bigcap\nolimits_{i\in[m-1]} H_i^+$. 
Then let us look at the following brief formula,

\begin{equation}\label{eq:Lawrence2}
\vol (P)=\sum_{\substack{\vbf: \text{ a vertex of}\\P\cap\{f>0\}}} N_\vbf. 
\end{equation}

The formula looks almost the same as Lawrence's original one, but ignores any vertex placed on $H_m$. Recall the geometric interpretation that $N_{\vbf}$ is the volume of a projected simplex. Intuitively, $N_{\vbf}$  for $f(\vbf)=0$ seems to vanish. If this is true, we should be able to handle some situations violating Lawrence's condition, i.e., if a non-simple vertex or a parallel edge  in $P$ is placed on an auxiliary plane, it is enough to simply ignore $N_{\vbf}$  as in (\ref{eq:Lawrence2}). We introduce a term to describe the situation as follows.
\begin{definition}
	A vertex or an edge of $P$ is \emph{degenerate} into a hyperplane $H_m$ if the vertex or the edge belongs to $H_m$, or  otherwise \emph{non-degenerate}. We usually omit indicating $H_m$ unless it will cause confusion.
\end{definition}

Now we verify the validity of the above intuition in the following theorem.    
\begin{theorem}\label{thm:mainLawrence-degenerate}
	Let a convex polyhedron $P =  \bigcap\nolimits_{i\in[m-1]} H_i^+\bigcap H_{m}^+$ where Lawrence's conditions hold only for non-degenerate vertices and edges.
	Then the volume is given by 
	\begin{equation}\label{eq:dege-form}
	\vol (P )
	=\mathlarger\sum_{\substack{I\subset[m-1]\\|I|=n}}
	~~\mathlarger\sum_{\vbf\in H_I}
	\frac{(-1)^{ \frac{n(n+1)}{2}} (g_m(\vbf) {{A}}_{[n]}^{I})^n}{n!| {{A}}_{[n]}^{I}| \prod\limits_{t\in I}  {{A}}_{[n]}^{I{\setminus t \cup m}}}.
	\end{equation}
	
\end{theorem}	

\begin{remark}
	Note that the formula itself is exactly the same as Theorem \ref{thm:mainLawrence} except the term $N_\vbf$ for any non-degenerate $\vbf$ isn't contained in the summation. 
	The number of hyperplanes involved in the computation is also reduced by one.
	In particular,
	the two requirements are practically weakened, i.e., it is enough to check whether simplicity and non-parallelness are fulfilled only at non-degenerate vertices and edges.
	
\end{remark}

\begin{proof}
	The proof consists of two steps. We first prove it for the case violating only  non-parallelness, i.e. we assume that $P$ has only simple vertices. 
	Consider a vector $\qbf$ that is not parallel to any edge of $P$ and choose an auxiliary function $f_\epsilon$ as follows,
	\begin{equation}\label{eq:f_epsilon}
	f_\epsilon(\xbf):= (\abf_m + \epsilon \qbf)\cdot \xbf + r_m.
	\end{equation} 
	Then we can apply Theorem \ref{thm:mainLawrence} to obtain $\vol(P)$ by putting 
	\[H_{m+1}:=\{\xbf \mid f_\epsilon(\xbf)=0\}\] as an auxiliary plane. 
	Let the volume of $P$ and the value of Theorem \ref{thm:mainLawrence} 
	be denoted by $\vol(P)$ and  $\vol_\epsilon(P)$ respectively. 
	By Lawrence, we have 
	$
	\vol(P)=\vol_\epsilon(P) \text{ for any } \epsilon>0
	$ 
	and the $\vol_\epsilon(P)$ is exactly represented by
	
	$$\begin{aligned}
	\vol_\epsilon (P )
	&=\mathlarger\sum_{\substack{I\subset[m]\\|I|=n}}
	~~\mathlarger\sum_{\vbf\in H_I}
	\frac{(-1)^{ \frac{n(n+1)}{2}} (g_{m+1}(\vbf) {{A}}_{[n]}^{I})^n}{n!| {{A}}_{[n]}^{I}| \prod\limits_{t\in I}  {{A}}_{[n]}^{I{\setminus t \cup (m+1)}}}.
	\end{aligned}
	$$
	One needs to be careful here: the number of hyperplanes increases by one from the count of Theorem  \ref{thm:mainLawrence}, i.e. $m\mapsto m+1$.
	We decompose the summation into two parts according to whether $\vbf$ belongs to $H_{m}$ or not. 
	\begin{align}
	=\mathlarger\sum_{\substack{I\subset[m-1]\\|I|=n}}
	~~\mathlarger\sum_{\vbf\in H_I}
	\frac{(-1)^{ \frac{n(n+1)}{2}} (g_{m+1}(\vbf) {{A}}_{[n]}^{I})^n}{n!| {{A}}_{[n]}^{I}| \prod\limits_{t\in I}  {{A}}_{[n]}^{I{\setminus t \cup (m+1)}}}\label{eq:firsthalf}\\
	+\mathlarger\sum_{\substack{\vbf\in H_{I\cup m}\\|I|=n-1}}
	\frac{(-1)^{ \frac{n(n+1)}{2}} (g_{m+1}(\vbf) {{A}}_{[n]}^{I\cup m})^n}{n!| {{A}}_{[n]}^{I\cup m}| \prod\limits_{t\in I\cup m}  {{A}}_{[n]}^{I\cup m{\setminus t \cup (m+1)}}}\label{eq:secondhalf}
	\end{align}
	
	Let us look at the first part  (\ref{eq:firsthalf}). The summand $N_\vbf^\epsilon$ for each vertex $\vbf\in H_I\setminus H_m$ goes to $N_\vbf$ in (\ref{eq:dege-form}), i.e.,
	\begin{align*}
	\lim_{\epsilon\rightarrow 0} N_\vbf^\epsilon =&
	\lim_{\epsilon\rightarrow 0}\mathlarger\sum_{\substack{I\subset[m-1]\\|I|=n}}
	~~\mathlarger\sum_{\vbf\in H_I}
	\frac{(-1)^{ \frac{n(n+1)}{2}} (g_{m+1}(\vbf) {{A}}_{[n]}^{I})^n}{n!| {{A}}_{[n]}^{I}| \prod\limits_{t\in I}  {{A}}_{[n]}^{I{\setminus t \cup (m+1)}}}\\
	=&\mathlarger\sum_{\substack{I\subset[m-1]\\|I|=n}}
	~~\mathlarger\sum_{\vbf\in H_I}
	\frac{(-1)^{ \frac{n(n+1)}{2}} (g_{m}(\vbf) {{A}}_{[n]}^{I})^n}{n!| {{A}}_{[n]}^{I}| \prod\limits_{t\in I}  {{A}}_{[n]}^{I{\setminus t \cup m}}}=N_\vbf. 
	\end{align*}
	Therefore the first part (\ref{eq:firsthalf}) goes exactly to the formula of (\ref{eq:dege-form}) as $\epsilon \to 0$. 
	
	Now we show that the second part (\ref{eq:secondhalf})  goes to zero as $\epsilon \to 0$. 
	For $\vbf\in H_m$, we have $g_{m+1}(\vbf)=f_\epsilon(\vbf)=\epsilon \qbf\cdot \vbf$ and  
	\[
	{{A}}_{[n]}^{I\cup m{\setminus t \cup (m+1)}}=\det(\abf_{i_1},\ldots,\abf_{i_{t-1}},\abf_{i_{t+1}},\ldots,\abf_m,\abf_m+ \epsilon \qbf).
	\] 
	From this, we obtain	
	\begin{align*}
	\prod\limits_{t\in I\cup m} {{A}}_{[n]}^{I\cup m{\setminus t \cup (m+1)}}
	=& \epsilon^{n-1}{{A}}_{[n]}^{I\cup m}\prod_{t=1}^{n-1}\det(\abf_{i_1},\ldots,\abf_{i_{t-1}},\abf_{i_{t+1}},\ldots,\abf_m,\qbf)\\
	&+\epsilon^{n}\det(\abf_{i_1},\ldots,\abf_{i_{n-1}},\qbf)\prod_{t=1}^{n-1}\det(\abf_{i_1},\ldots,\abf_{i_{t-1}},\abf_{i_{t+1}},\ldots,\abf_m,\qbf).
	\end{align*}
Using these computations, we know that the denominator  and 
	the numerator of $N_\vbf^\epsilon$ in the expression (\ref{eq:secondhalf}) have  dominating terms of $\epsilon^{n-1}$ and $\epsilon^{n}$, respectively. Therefore we conclude that the second part goes to zero as $\epsilon\to 0$.

	Now we can freely use the formula (\ref{eq:dege-form})
	for a simple polytope $P$ even cases violating the non-parallel condition as long as the edge is degenerate into $H_m$.   
	The remaining case is  that $P$ has a  non-simple vertex which is degenerate into  $H_m$. 
	Let us consider a truncated polytope $P^\epsilon$ obtained from $P$ by cutting out a sufficiently small neighborhood of each vertex $\vbf$ by an additional hyperplane $H_\vbf^\epsilon$. More precisely, 
	a non-simple vertex $\vbf$ is removed from $P$ and a new facet $H^\epsilon_\vbf$ appears with several simple vertices $\wbf$'s such that  
	\begin{inparaenum}[(i)]
		\item the normal vector of $H_\vbf^\epsilon$ is fixed as  $\epsilon \to 0$,
		\item  when  $\epsilon =0$, 
		$H_\vbf^\epsilon$ contains $\vbf$, and 
		\item $f(\wbf) < \epsilon$ for any new vertex $\wbf$ from the truncation by $H^\epsilon_\vbf$. 	
	\end{inparaenum}
	
	Let $\vol^f(P)$ be the value of the formula (\ref{eq:dege-form}). 
	We estimate $\vol^f(P)$ as follows,
	\[
	\begin{aligned}
	\vol(P)-\vol^f(P) &= \vol(P)-\vol^f(P^\epsilon) + \sum_\wbf N^\epsilon_\wbf\\
	\end{aligned}
	\]
	where $N^\epsilon_\wbf$ is the summand of $\vol^f(P^\epsilon)$ at $\wbf$ of (\ref{eq:dege-form}).
	Note that we know that $\vol^f(P^\epsilon)=\vol(P^\epsilon)$ since $P^\epsilon$ fulfills the simplicity condition. 
	So the difference  $\vol(P) - \vol^f(P^\epsilon)$ becomes
	$\vol(P) - \vol(P^\epsilon)$ which obviously goes to zero as $\epsilon \to 0$ since they are geometric Euclidean volumes. 
	Moreover,  $N^\epsilon_\wbf$ also goes to zero since $N^\epsilon_\wbf= c \epsilon^n $ for some constant $c$ because the normal vector of $H_\vbf^\epsilon$ is fixed while $\epsilon \to 0$.    
	Therefore we conclude that even if $P$ contains a non-simple vertex,  $\vol(P)=\vol^f(P)$ because of the degeneracy of the vertex $\vbf$.       
\end{proof}

{
	\subsection{$\epsilon$-pertubation for general polytopes}\label{sec:epsilonpert}
	Sometimes we want to consider the volume of a general polytope that might violate  Lawrence's conditions.	
	For example, a way to apply Lawrence's formula to a non-simple polytope is to decompose  a non-simple vertex into projected simplices along ``lexicographic rule" (\cite{gritzmann_complexity_1994}, \cite{lawrence_polytope_1991}).  
	
	In the previous section, our formula weakened  Lawrence's conditions in a specific situation. However, the strategy using a limiting process can be, in principle, utilized in general.  
	Here we discuss such a ubiquitous method, called \emph{$\epsilon$-perturbation},  for  cases  violating either simplicity or non-parallelness. Let us say a polytope is \emph{good} if it satisfies Lawrence's conditions, and otherwise \emph{bad}.
	
	The following two facts constitute the essential justification for the  $\epsilon$-perturbation method. 
	
	\begin{enumerate}
		\item The volume of a clipped hypercube is a continuous function on the parameter spaces of the  hyperplanes. 
		\item Lawrence's condition is an open condition. Precisely, the set of bad  polytopes is a Zariski-closed subset in the parameter space.
	\end{enumerate}
	Therefore 
	we can always find a good polytope which has only an arbitrarily small  difference from the original bad polytope. Moreover the continuity of volume tells us that the volume difference is also arbitrarily small. 
	
	Let us explain the strategy more precisely as follows. 
	For a bad polytope $P$,  
	let us add a perturbation variable $\epsilon$ (or sometimes more hyerperplanes) into a system of hyperplanes such that the perturbed polytopes  $P_\epsilon$ are  always good for any sufficiently small $\epsilon >0 $. Then we can obtain the volume of $P$ as the limit as $\epsilon \rightarrow 0$ by only using the formula  for good polytopes.
	
	\begin{remark}
		In fact, the strategy itself is nothing special but it is usually not easy because the analysis of a precise limit is difficult in many cases. In particular, to apply this, a sufficiently concrete expression is required as in Theorem \ref{thm:mainLawrence-degenerate} and Section \ref{sec:concereteexample}.
		We are going to give more examples of $\epsilon$-perturbation  in Section \ref{sec:combidentity} and the Appendix. 
		In a follow-up article \cite{cho_exact_????},  a concrete exposition of $\epsilon$-perturbation for several cases are investigated.  
	\end{remark}

}
\section{Volume formulas for a clipped hypercube}\label{sec:mainproof}   

Now we derive volume formulas for a clipped hypercube. A clipped hypercube has two types of hyperplanes: the facets of $[0,1]^n$ and the clipping hyperplanes $H_1,\dots, H_m$. We use calligraphic notation $\Hcal, \Acal$ and $\Ical$ when we consider both types of hyperplanes in a unified way.

\subsection{Hyperplanes and indices for clipped hypercubes}
\label{sec:clippedcubedefinition}
An $n$-dimensional unit hypercube $[0,1]^n$ is given by $2n$  half spaces,
\begin{align*}
x_i\geq 0  ~~\text{ and }~~ x_i\leq 1  ~~\text{ for }~~ i \in [n].
\end{align*}
Let $P$ be a hypercube clipped by hyperplanes $H_1, \dots H_m$, i.e. 
\begin{align*}
P := [0,1]^n \cap \bigcap\limits_{i=1}^m H_i^+.
\end{align*} 
From now on, we use separate  notation $\Hcal_i$ and $H_i$ for  hyperplanes. The former $\Hcal_i$ indicates a hyperplane in the set of all hyperplanes of $P$, while the later $H_i$ indicates a hyperplane in the subset of only the clipping hyperplanes.
We follow the same notational conventions for $\Hcal$, $\Hcal^+$ and $\Hcal_i$ and so on, as in Section \ref{sec:hyperplaneNotation}.
Let us consider the set of half spaces  $\Hcal_i^+$. In order to apply Theorem \ref{thm:mainLawrence}, consider
\begin{align}
\Hcal_i^+ = \left\{ 
\begin{aligned}
\{ \xbf \in \reals^n ~&|~  \ebf_{(i+1)/2}\cdot \xbf \geq 0 \}, &\text{ if }& i \in [2n] \text{ is odd}\\
\{ \xbf \in \reals^n ~&|~ (-\ebf_{i/2})\cdot \xbf +1 \geq 0 \}, &\text{ if }& i \in [2n] \text{ is even}\\ 
&H_{i-2n}^+ ,&\text{ if }& i \in [2n+m]\setminus[2n]. \\
\end{aligned}	
\right.
\end{align}
We also consider a big coefficient matrix $\Acal$  and the index set $\Ical=[2n+m]$ for $\Hcal_i$, where the matrix $A$ is given by hyperplanes $H_i$ and the index set $I=[m]$.
We construct an $n \times (2n+m)$ matrix $\Acal$ which is determined by $\Hcal$ as follows. Let us keep in mind the canonical embedding $I$ into $\Ical$ by $i\mapsto i+2n$.
\begin{align}
\Acal 
&=(~\ebf_1 , -\ebf_1 , \ebf_2 , -\ebf_2 , \dots , \ebf_n , -\ebf_n ~|~ A ~ ) \notag \\
&=\left[\begin{array}{ccccccccccc}
1&-1&&&&&&a_{11}&&a_{1,m-1}&a_{1,m}\\
&&1&-1&&&&a_{21}&&a_{2,m-1}&a_{2,m}\\
&&&&\ddots&&&\vdots&\cdots&\vdots&\vdots\\
&&&&&1&-1&a_{n1}&&a_{n,m-1}&a_{n,m}\\
\end{array}\right]
\end{align}

For a simple vertex $\vbf$, there is  an index set $\Ical$ tracking the $n$  hyperplanes $\Hcal_{i_1},\Hcal_{i_2},\dots,\Hcal_{i_n}$ making up  the vertex $\vbf$, i.e. 
\begin{align*}
\vbf=\bigcap_{k\in \Ical}\Hcal_k ~\text{ for }~ \Ical =\Ical(\vbf)=\{i_1,\dots,i_n\} \subset [2n+m],
\end{align*}
where  $\Ical$ is well-ordered, $ i_1 < i_2 < \dots < i_n$.
We decompose $\Ical$ into two parts,   
\begin{align}\label{eq:I01Istar}
\Ical = \Ical_{01} \cupdot \Ical_*
\end{align}
with 
$\Ical_{01} = \{i\in \Ical | i\in [2n] \}$ and $\Ical_{*} = \{i\in \Ical | i\in [2n+m]\setminus[2n] \}$.
We can check that
\begin{align}\label{eqn:AI*}
(\Acal)^{\Ical_*}_{[n]} = (A)^I_{[n]} .
\end{align}
The following lemma says that each vertex has a natural grading from $F^d$. Our volume formula is a summation over the grading.
\begin{lemma} \label{lem:vstar_index}
	For a simple vertex $\vbf$ of a clipped hypercube $P$, there is an index set $I \subset [n]$ such that $\vbf \in F^d \cap H_I$. Moreover,	
	$
	|\starv| = |I|=d
	$.
\end{lemma}
\begin{proof}
	Let  $\vbf = (v_1,v_2,\dots,v_n) \in F^d$, then $|\bulletv|=|\Ical_{01}|=n-d$ because $v_i = 0 \text{ or } 1$ if and only if $\vbf$ intersects a hyperplane of form $x_i = 0$ or $x_i = 1$.
	Since $ \Ical_* = I + 2n$ and $| \Ical_* | = |I|$, the lemma is obvious by (\ref{def:v01}) and (\ref{eq:I01Istar}).
\end{proof}
Let $\vbf \in F^d$ and let us consider this in  more detail,
\begin{align}
\vbf \in \bigcap_{i \in \zerov} \Hcal_{2i-1} \cap \bigcap_{i \in \onev} \Hcal_{2i}.
\end{align}	
Let $\begin{aligned}
\Ical_{0} := \{i \in \Ical_{01} ~|~ i \text{ is odd }\} \text{ and }~ \Ical_{1} := \{
i \in \Ical_{01} ~|~ i \text{ is even }\} 
\end{aligned}$. Then it  immediately follows that 
\begin{align}\label{deff}
\Ical_0  ={  (2 \cdot \zerov - 1 ) } \text{ and }~
\Ical_1 ={   2 \cdot \onev. } 
\end{align}
Lastly, we also obtain the following  lemma. 
\begin{lemma}\label{lem:Av01I01}
	$(\Acal)^{\Ical_{01}}_{\bulletv}$ is a diagonal matrix with
	\begin{align*}
	((\Acal)^{\Ical_{01}}_{\bulletv})_{i,i} = \left\{ 
	\begin{aligned}
	1 & ~~\textup{ if } ~ v_{\bullet_i(\vbf)} = 0 \\
	-1 & ~~\textup{ if }  ~ v_{\bullet_i(\vbf)} = 1. 
	\end{aligned}
	\right.
	\end{align*}	
	In particular, $\Acal_{\bulletv}^{\Ical_{01}} = (-1)^{|\onev|}$.	 
\end{lemma}

\begin{proof}
	The hyperplanes $\Hcal_{2i}$ and $\Hcal_{2i-1}$ are parallel and never intersect. Hence, for any $t \in [n]$,  ${2t}$ and ${2t-1}$ cannot be contained in ${\Ical_{01}}$ at the same time. By the definition of $\Acal$ and (\ref{deff}), the matrix is diagonal with $|\zerov|$ entries equal to $1$ and $|\onev|$ entries equal to  $-1$.  	
\end{proof}

\subsection{Good clipping conditions}\label{sec:goodclipping}
We consider two more explicit conditions called good clipping conditions in which a volume formula is applicable. They look almost identical to Lawrence's two conditions in a concrete form, but are not exactly the same.
First of all,  we don't need any requirement for a degenerate vertex or  edge since Theorem \ref{thm:mainLawrence-degenerate} doesn't require them. So a polytope violating simplicity or non-parallelness may satisfy our good clipping conditions.
Not only that, but Lawrence's condition doesn't imply the good clipping conditions. To see this concretely, see Example \ref{exam:simplebutnotgood}.

\begin{proposition}[Good clipping conditions]\label{prop:goodclippingass}
	If a clipped hypercube $P$ 
	satisfies the following assumption :
	\begin{enumerate}[(A)]
		\item For any  $I \subset [m-1]$, $F^{|I|-1} \cap H_{I}=\varnothing$,
		\item  For any $I\subset[m-1] $ and $\vbf \in F^{|I|}\cap H_I$, $\prod\limits_{t\in I} A_{\starv}^{I{\cup m\setminus t}} \prod\limits_{t \in \bulletv} A_{\starv \cup t} ^{I{\cup m}}\neq 0$,
	\end{enumerate}
	then  Lawrence's two conditions hold for non-degenerate vertices.  
\end{proposition}

\begin{proof}
	
	Let us see that (A) implies that every non-degenerate vertex is simple.
	If a non-simple vertex $\vbf$ exists then there are at least $n+1$   hyperplanes intersecting   $\vbf$, which  meet $k$  hyperplanes of $[0,1]^n$ and $H_{i_1},\dots ,H_{i_{n-k+1}}$ for $0\leq k \leq n+1$. Then  $\vbf \in F^{n-k} \cap H_I $ with $I=\{i_1,\dots,i_{n-k+1}\}$. This conflicts with (A).
	Condition  (B) means the volume $N_\vbf$ of the projected simplex at $\vbf$ is finite, which is equivalent to Lawrence's non-parallelness. 
\end{proof}

Let us give an example of a polytope $P$ which violates the good clipping condition (A) but  satisfies Lawrence's simplicity. Note  the non-parallel condition always guarantees good clipping condition (B).

\begin{example}\label{exam:simplebutnotgood}
	There is  a clipped hypercube which is simple but violates the good clipping condition (A). Consider $P=[0,1]^3 \cap H_1^+ \cap H_2^+$ where
	\begin{align*}
	H_1^+:&=\{ (x_1,x_2,x_3) \in \reals^3| x_1+x_2+2 x_3 \geq 1\}, \\
	H_2^+:&=\{ (x_1,x_2,x_3) \in \reals^3| x_1+x_2 \leq 1  \}.
	\end{align*}
	We can see that  all vertices are  simple, but  
	\begin{align*}
	H_1 \cap F^{1-1} &= \{(1,0,0),(0,1,0)\}\neq \varnothing, \\
	H_2 \cap F^{1-1} &= \{(1,0,0),(0,1,0),(1,0,1),(0,1,1)\}\neq \varnothing. 
	\end{align*}
	In this example, we observe that there can be a simple vertex in a clipped hypercube, which lies in the intersection of more than $n+1$  hyperplanes in $\Hcal$. When we enumerate  hyperplanes at  a vertex of a clipped hypercube, the hyperplanes in $[0,1]^n$ are chosen by a default, but some of them are sometimes  unnecessary because of clipping  hyperplanes in $H$.
\end{example}
Finally, we remark that Conditions (A) and (B) are also open conditions  like Lawrence's conditions.

\subsection{Several  volume formulas}\label{sec:severalversion}
We now present our main formulas. Refer to Section \ref{sec:notation}  and  \ref{sec:clippedcubedefinition} for notation.
We will postpone the proof until later. 
If one wants to obtain each summand  up to sign, the evaluation is much simpler. However, we have to keep track of correct signs to derive an explicit formula. The practical difficulty of the proof comes from handling the complicated signs and hence we need to prepare something to deal with the signs properly. 

Here, we are going to discuss several different presentations of the same formula and the relevant lemmas to convert among the formulas using different conventions.
Let us look at a first expression for the volume of clipped hypercube as follows.
\begin{theorem}\label{thm:maincup}
	The volume of a hypercube clipped by $m$  halfspaces $H_1,H_2,\dots, H_m$  
	satisfying good clipping conditions is given by 
	\begin{align*}
	\vol ( [0,1]^n \cap {H^+} )
	= \sum_{ I\subset[m-1]}
	\;\;\sum_{ \vbf \in F^{|I|} \cap H_{I}}
	\frac{(-1)^{|\zerov|+ \| \starv\| } (g_m(\vbf)A_{\starv}^{I})^n} {n!|A_{\starv}^{I}| \prod\limits_{t\in I} A_{\starv}^{I{\cup m\setminus t}} \prod\limits_{t \in \bulletv} A_{\starv \cup t} ^{I \cup m } }.
	\end{align*}		
	%
\end{theorem}
We consider another expression which uses joining union  $\vee$ instead of ordered union $\cup$ as follows. The only difference between them is the sign computation for $\frac{|I|(|I|+1)}{2}$ and $\| \starv \|$. 
\begin{theorem}\label{thm:mainvee}
	
	The volume of a hypercube clipped by $m$  halfspaces $H_1,H_2,\dots, H_m$  
	satisfying good clipping conditions is given by 
	\begin{align*}
	\vol ( [0,1]^n \cap {H^+} )
	= \sum_{ I\subset[m-1]}
	\;\;\sum_{ \vbf \in F^{|I|} \cap H_{I}}
	\frac{(-1)^{|\zerov|+ \frac{|I|(|I|+1)}{2}}  (g_m(\vbf)A_{\starv}^{I})^n} {n!|A_{\starv}^{I}| \prod\limits_{t\in I} A_{\starv}^{I{\vee m\setminus t}} \prod\limits_{t \in \bulletv} A_{\starv \vee t} ^{I \vee m } }.
	\end{align*}		
\end{theorem}

The following lemma implies that the above two theorems are equivalent. 
\begin{lemma}\label{lem:veecup}
	$$
	\prod_{t\in \bulletv} A_{\bulletv \cup t}^{I{\cup m}}
	=(-1)^{\|\starv\|- \frac{|I|(|I|+1)}{2}}\prod_{t \in \bulletv} 
	A_{\starv \vee t}^{I{\vee m}}.
	$$
\end{lemma}
\begin{proof}
	The set $\bulletv = [n]\setminus \starv$ is divided into $|I|+1$  (possibly empty) blocks,  
	\begin{align*}
	&\{1,2,\ldots,i_1-1\}, \\
	&\{i_1+1,i_1+2,\ldots,i_2-1\},\\ 
	&~~~~~~~~~~~\vdots \\
	&\{i_{|I|}+1,i_{|I|}+2,\ldots,n\}. 
	\end{align*}
	Whenever $t \in \bulletv$ is placed in each block, it requires $|I|,|I|-1,\ldots,1,0$  transpositions. Hence the total number of transpositions is 
	\begin{align*}
	&(i_1-1)|I|+(i_2-i_1-1)(|I|-1)+\cdots+(i_{|I|}-i_{|I|-1}-1)1+(n-i_{|I|})0\\ &=i_1+i_2\cdots+i_{|I|}-(|I|+(|I|-1)+\cdots+1)\\
	&=\|\starv\|-\frac{|I|(|I|+1)}{2}. 
	\end{align*}
\end{proof}

Recall that the ordered union is commutative but the joining union is not, i.e. 
$ I \cup m = m \cup I $ but $I \vee m \neq m \vee I$.
Thus, in  Theorem \ref{thm:maincup} the  expression order of the union operations doesn't matter. But if one takes the joining union as in Theorem \ref{thm:mainvee}, there are several choices of expressions  because it is sensitive to changing order. 
The following lemma
shows that  orders of the expressions  affects a few things in Theorem \ref{thm:mainvee}.
\begin{lemma}\label{lem:veecommu}
	$$	\prod\limits_{t\in I} A_{\starv}^{I{\vee m\setminus t}} =
	\prod\limits_{t\in I} A_{\starv}^{m \vee  I \setminus t} = \prod\limits_{t\in I} A_{\starv}^{I{\cup m\setminus t}} =
	\prod\limits_{t\in I} A_{\starv}^{m \cup  I \setminus t},
	$$	
	$$
	\prod\limits_{t\in \bulletv} A_{\starv\vee t }^{  I \vee m}
	=\prod\limits_{t\in \bulletv} A_{t \vee \starv}^{ m \vee I}	
	= (-1)^{|I|(n-|I|)} \prod\limits_{t\in \bulletv} A_{\starv\vee t }^{ m \vee I}
	= (-1)^{|I|(n-|I|)}  \prod\limits_{t\in \bulletv} A_{ t \vee \starv }^{  I \vee m}.	
	$$	
\end{lemma}
\begin{proof}
	This is obvious by	Lemma \ref{lem:vstar_index} along with 
	\begin{align*}
	\sigma (m \cup I, I \cup m) &= 1 ~~\text{ and }~~ \sigma( m \vee I,  I \vee m)= (-1)^{|I|}. 
	\end{align*}
\end{proof}

Lastly, we write down another version of formula that uses the joining union. The effect of making another choice of order can be computed easily by Lemma \ref{lem:veecommu}. Here we use $m\vee I$ instead of $I \vee m$ in Theorem \ref{thm:mainvee} and the effect is summarized  via swapping $\frac{|I|(|I|+1)}{2}$ for $\frac{|I|(|I|-1)}{2} + n |I|$.
\begin{theorem}\label{thm:mainvee2}
	Under the same 	hypotheses as Theorem \ref{thm:mainvee},
	\begin{align*}
	\vol ( [0,1]^n \cap {H^+} )
	= \sum_{ I\subset[m-1]}
	\;\;\sum_{ \vbf \in F^{|I|} \cap H_{I}}
	\frac{(-1)^{|\zerov|+ \frac{|I|(|I|-1)}{2}+n|I|}  (g_m(\vbf)A_{\starv}^{I})^n} {n!|A_{\starv}^{I}| \prod\limits_{t\in I} A_{\starv}^{I{\vee m\setminus t}} \prod\limits_{t \in \bulletv} A_{\starv \vee t} ^{m \vee I} }.
	\end{align*}		
\end{theorem}

\subsection{Proof of the volume formula for clipped hypercubes }

In this section, we prove  Theorem \ref{thm:mainvee}. The other formulas of Theorem \ref{thm:maincup} and Theorem \ref{thm:mainvee2} can be derived from it as discussed in the previous section. 
For effective bookkeeping of complicated sign permutations, we introduce  the \emph{separating parity} $\Delta(I,J)$ of two indices $I \supset J$ as follows,
$$\Delta(I,J) := \sigma ( I ,(I \setminus J )\vee J ). $$

From the definition, we have the following lemmas about $\Delta$.

\begin{lemma}
	For any $I \subseteq [n]$,
	$$\Delta([n],I) = (-1)^{n|I|- \| I \| - \frac{|I|(|I|-1)}{2}}.$$
In particular, for $t \in [n]$
	$$\Delta([n],t)  = (-1)^{n-t}.$$ 
\end{lemma}
\begin{proof}
	Let $I=\{i_1,\dots,i_{|I|}\}$, then  
	\begin{align*}
	[n]\setminus I  &= \{ 1, \dots, i_1-1,\widehat{i_1}, i_1+1, \dots, i_{|I|}-1,  \widehat{i_{|I|}}, i_{|I|}+1, \dots,n \},  \\
	([n] \setminus I )\vee I & = \{ 1, \dots, i_1-1,  i_1+1\dots, i_{|I|}-1, i_{|I|}+1, \dots,n \} \vee \{{i_1},\dots ,{i_{|I|}} \}.
	\end{align*}
	We just count the number of transpositions.
	In order to shift each $i_t$ in $([n] \setminus I) \vee I$ into its original position in $[n]$, it requires $n - i_t - (|I|-t)$  transpositions.
	\begin{align*}
	\sum_{t\in [|I|]}  n - i_t - (|I|-t) &= n |I| - \|I\| -\frac{|I|(|I|-1)}{2}. 
	\end{align*}
\end{proof}

\begin{lemma}
	For any $I \subseteq [n]$,
	$$
	\prod_{i\in I} \Delta(I,i) = (-1)^{\frac{|I|(|I|-1)}{2}}.
	$$
\end{lemma}
\begin{proof}
	Let $I=\{i_1,i_1,\dots,i_{|I|}\}$. Each $i_t = i$ requires $|I|-t$  transpositions. Hence,	
	\begin{align*}
	(-1)^{|I|-1}(-1)^{|I|-2}\cdots(-1)^{|I|-|I|} &= (-1)^{\frac{|I|(|I|-1)}{2}}. 
	\end{align*}
\end{proof}

%
At the first step of the proof, the summation over $\Ical$ and $\Hcal_\Ical$ applied to Theorem \ref{thm:mainLawrence} is converted into a summation over $I$ and $H_I$ as follows.
\begin{align*}
\sum_{\substack{\Ical \subset [2n+m-1]\\|\Ical|=n}} ~~ \sum_{\vbf \in \Hcal_\Ical}
&= \sum_{\substack{\Ical_{01} \vee \Ical_*  \subset [2n+m-1]\\|\Ical_{01} \vee \Ical_*|=n}} ~~~ \sum_{\vbf \in \Hcal_{\Ical_{01}} \cap \Hcal_{\Ical_* }} \\
&= \sum_{\Ical_* \subset [2n+m-1]\setminus[2n]} 
~~~\sum_{\substack{\Ical_{01} \subset [2n]\\|\Ical_{01}|+|\Ical_*| = n }} ~~~ \sum_{\vbf \in \Hcal_{\Ical_{01}} \cap \Hcal_{\Ical_* }} \\
&= \sum_{\Ical_* \subset [m-1]+2n } 
~~~\sum_{\vbf \in F^{|\Ical_{*}|} \cap \Hcal_{\Ical_* }} \\
&= \sum_{I \subset [m-1] } 
~~~\sum_{\vbf \in F^{|I|} \cap H_I } 
\end{align*}

For the sake of convenience, let $\Delta$ denote the separating parity of  $[n]$ and $\starv$, i.e.,
\begin{align}\label{eqn:seperatingv}
\Delta := \Delta([n],\starv) 
&=\sigma(\bulletv \cup \starv , \bulletv \vee \starv ).  
\end{align}
Then we derive several relations between minors of $A$ and  $\Acal$.

\begin{proposition}\label{prop:AtoA}
	$$\Acal^\Ical_{[n]} ~=~  (-1)^{|\onev|} ~ \Delta  ~ A^I_{\starv} $$
\end{proposition}
\begin{proof}
	By Lemma \ref{lem:Av01I01} and (\ref{eqn:AI*}),
	\begin{align*}
	\Acal^\Ical_{[n]} ~&{=}~~ \Acal ^{\Ical_{01}\vee \Ical_*} _{\bulletv \cup *_{\vbf} } ~{=}~ \Delta ~ \Acal ^{\Ical_{01}\vee \Ical_*} _{\bulletv \vee \starv }\\
	~&=~~ \Delta ~\Acal ^{\Ical_{01}}_{\bulletv} ~ \Acal^{ \Ical_*}_{ \starv } 
	~=~~  (-1)^{|\onev|} ~ \Delta  ~ A^I_{\starv}. 
	\end{align*}	
\end{proof}
Recall  Lemma \ref{lem:Av01I01} and (\ref{eqn:AI*}); the matrix $\Acal_{[n]}^{\Ical_{01}}$ is diagonal and  $\Acal_{[n]}^{\Ical_{*}} = A_{[n]}^I$.  Hence we need to decompose  $\Acal^{\Ical \setminus t \cup m}_{[n]}$ into  $\Ical_{01}$ and $\Ical_*$ as follows.
\begin{align}\label{eqn:a_decompose}
\prod_{t\in \Ical}\Acal^{\Ical \setminus t \cup m}_{[n]} ~&=~ \prod_{t\in \Ical_{01}}\Acal^{(\Ical_{01} \setminus t) \cup \Ical_* \cup m}_{[n]} ~ \prod_{t\in \Ical_*}\Acal^{\Ical_{01}  \cup ( \Ical_* \setminus t)  \cup m}_{[n]}.
\end{align}
The case of $\Ical_* \setminus t$ is proved  in the same way.
\begin{proposition} \label{prop:A_I*toA}
	\begin{align*}
	\prod_{t\in \Ical_*} \Acal^{\Ical_{01} \vee ( \Ical_* \setminus t \vee m )}_{[n]} ~&=~  
	(-1)^{|I||\onev|} ~  \Delta^{|I|} ~\prod_{t\in I }  A^{I\setminus t \vee m}_{\starv}.
	\end{align*}	
\end{proposition}
The case  $\Ical_{01} \setminus t$ is much more complicated than the previous case. 
\begin{proposition}\label{prop:A_01toA}
	\begin{align*}
	\prod_{t\in \Ical_{01}} ~
	\Acal^{(\Ical_{01} \setminus t) \vee \Ical_* \vee m}_{[n]} 
	~
	&=~ (-1)^{|\starv||\zerov|+|\onev||\zerov| +|\starv||\onev|+ \frac{|\bulletv|(|\bulletv|-1)}{2} }
	\Delta^{|\bulletv|}
	\prod_{t\in \bulletv} 
	A^{I_* \vee m}_{ \starv \vee t }.
	\end{align*}	
\end{proposition}
\begin{proof}
	First, divide $\Ical_{01}$  into $\Ical_0$ and $\Ical_1$, which correspond to $\zerov$ and $\onev$ respectively.
	\begin{align*}
	\prod_{t\in \Ical_{01}} ~
	\Acal^{(\Ical_{01} \setminus t) \vee \Ical_* \vee m}_{[n]} 
	&=~ \prod_{t\in \Ical_{01}}
	~\Acal^{(\Ical_{01} \setminus t) \vee \Ical_* \vee m}_{{\bulletv \cup \starv }} \\
	~
	&=~\prod_{t\in \Ical_{01}} ~\Delta  ~ \Acal^{(\Ical_{01} \setminus t) \vee \Ical_* \vee m}_{{\bulletv \vee \starv }} \\
	~
	&=~\Delta^{|\Ical_{01}|} ~\prod_{t\in \Ical_{0}}   ~ \Acal^{(\Ical_{01} \setminus t) \vee \Ical_* \vee m}_{{\bulletv \vee \starv }}
	~\prod_{t\in \Ical_{1}}   ~ \Acal^{(\Ical_{01} \setminus t) \vee \Ical_* \vee m}_{{\bulletv \vee \starv }} \\
	~
	&=~\Delta^{|\Ical_{01}|} ~
	\prod_{t\in \zerov}   ~ \Acal^{(\Ical_{01} \setminus (2t-1)) \vee \Ical_* \vee m}_{{\bulletv \vee \starv }}
	~\prod_{t\in \onev}   ~ \Acal^{(\Ical_{01} \setminus 2t) \vee \Ical_* \vee m}_{{\bulletv \vee \starv }}.
	\end{align*}
	Each term is computed as follows.
	\begin{align}
	\prod_{t\in \zerov}   ~ \Acal^{(\Ical_{01} \setminus (2t-1)) \vee \Ical_* \vee m}_{{\bulletv \vee \starv }}
	~
	&=~
	\prod_{t\in \zerov} \Delta(\bulletv,t)   ~ \Acal^{(\Ical_{01} \setminus (2t-1)) \vee \Ical_* \vee m}_{{(\bulletv \setminus t) \vee t \vee \starv }} \notag\\
	~
	&=~
	(-1)^{|\starv||\zerov|}\prod_{t\in \zerov} \Delta(\bulletv,t)    \Acal^{(\Ical_{01} \setminus (2t-1)) \vee \Ical_* \vee m}_{{(\bulletv \setminus t) \vee \starv }\vee t }. \notag\\
	\intertext{
		By Lemma \ref{lem:Av01I01}, we have   
		$
		\Acal^{\Ical_{01}}_{{\bulletv }} = 
		\left\{
		\begin{aligned}
		\Acal^{\Ical_{01} \setminus (2t-1)}_{{\bulletv \setminus t }} ~~~&\text{ if $t \in \zerov$}\\
		-\Acal^{\Ical_{01} \setminus 2t}_{{\bulletv \setminus t }} ~~~&\text{ if $t \in \onev$}
		\end{aligned}
		\right.
		$	
		. Then, }
	&=~
	(-1)^{|\starv||\zerov|}\prod_{t\in \zerov} \Delta(\bulletv,t)    \Acal^{\Ical_{01} \setminus (2t-1)}_{{\bulletv \setminus t }}
	A^{\Ical_* \vee m}_{ \starv \vee t }\notag\\
	&=~
	(-1)^{|\starv||\zerov|}\prod_{t\in \zerov} \Delta(\bulletv,t)    
	(-1)^{|\onev|}
	\Acal^{\Ical_* \vee m}_{ \starv \vee t } \notag\\
	&=~
	(-1)^{(|\starv|+|\onev|)|\zerov|}\prod_{t\in \zerov} \Delta(\bulletv,t)    
	\Acal^{I_* \vee m}_{ \starv \vee t } \label{eqn:I01_2t-1}.
	\end{align}

	Similarly,
	\begin{align}
	\prod_{t\in \onev}   ~ \Acal^{(\Ical_{01} \setminus 2t) \vee \Ical_* \vee m}_{{\bulletv \vee \starv }}
	~
	&=~
	(-1)^{(|\starv|+|\onev|)|\onev| - |\onev|}\prod_{t\in \onev} \Delta(\bulletv,t)    
	A^{I_* \vee m}_{ \starv \vee t } \label{eqn:I01_2t}.
	\end{align}
	
	Take (\ref{eqn:I01_2t-1}) and (\ref{eqn:I01_2t}) together to complete the proof. 
	\begin{align*}
	\prod_{t\in \Ical_{01}} ~
	\Acal^{(\Ical_{01} \setminus t) \vee \Ical_* \vee m}_{[n]} 
	~
	&=~ (-1)^{ (|\starv|+|\onev|)(|\zerov|+|\onev|) -|\onev|}\Delta^{|\Ical_{01}|}
	\prod_{t\in \bulletv} \Delta(\bulletv,t)    
	A^{I_* \vee m}_{ \starv \vee t }\\
	~
	&=~ (-1)^{|\starv||\zerov|+|\onev||\zerov| +|\starv||\onev|+ \frac{|\bulletv|(|\bulletv|-1)}{2} }
	\Delta^{|\bulletv|}
	\prod_{t\in \bulletv} 
	A^{I_* \vee m}_{ \starv \vee t }.
	\end{align*}
\end{proof}

We put the three propositions \ref{prop:AtoA}, \ref{prop:A_I*toA}
and \ref{prop:A_01toA} together. 
\begin{align*}
&\sum_{\substack{\Ical\subset[m-1]\\|\Ical|=n}}
~~\sum_{\vbf\in \Hcal_\Ical}
\frac{(-1)^{\frac{n(n+1)}{2}} (g_m(\vbf) {{\Acal}}_{[n]}^{\Ical})^n}{n!| {{\Acal}}_{[n]}^{\Ical}| \prod\limits_{t\in \Ical}  {{\Acal}}_{[n]}^{\Ical{\setminus t \cup m}}}\\
&= 
\sum_{I \subset [m-1] } 
\sum_{\vbf \in F^{|I|} \cap H_I }
\frac{ (-1)^{\frac{n(n+1)}{2}+n|\onev|+|I||\onev|+|\starv||\zerov|+|\onev||\zerov| +|\starv||\onev|+ \frac{|\bulletv|(|\bulletv|-1)}{2}}  \Delta^n  (g_m(\vbf)  A^I_{\starv} )^n }
{n!|A^I_{\starv}| \Delta^{|I|+\bulletv} \prod\limits_{t\in I }  A^{I\setminus t \vee m}_{\starv}
	\prod\limits_{t\in \bulletv} 
	A^{I_* \vee m}_{ \starv \vee t } }.
\end{align*}
We calculate the parity expression,
\begin{align*}
&\frac{n(n+1)}{2}+n|\onev|+|I||\onev|+|\starv||\zerov|+|\onev||\zerov| +|\starv||\onev|+ \frac{|\bulletv|(|\bulletv|-1)}{2} \\
&~~~\underset{( \text{mod } 2 )}{\equiv}\frac{n(n+1)}{2}+ \frac{|\bulletv|(|\bulletv|+1)}{2} -|\bulletv| +n|\onev|+|\starv||\zerov|+|\onev||\zerov| \\
&~~~\underset{( \text{mod } 2 )}{\equiv} \frac{|\starv|(|\starv|+1)}{2} + |\starv||\bulletv|+|\bulletv|+(n+|\zerov| )|\onev| +|\starv||\zerov| \\
&~~~\underset{( \text{mod } 2 )}{\equiv} \frac{|\starv|(|\starv|+1)}{2} + |\starv||\onev|+|\bulletv|+(|\starv|+|\onev| )|\onev| \\
&~~~\underset{( \text{mod } 2 )}{\equiv} \frac{|\starv|(|\starv|+1)}{2} + |\zerov|.
\end{align*}
This completes the proof of Theorem \ref{thm:mainvee}.

%
%
%
%

\section{More explicit formulas for $m \leq 3$ }\label{sec:concereteexample}
For the case of a small number of hyperplanes, we  derive more concrete formulas in a fully explicit way. We expect for these formulas to be more accessible for a broader  range of readers. Furthermore  such an elementary and precise formulations play a crucial role in obtaining combinatorial identities in Section \ref{sec:combidentity}.

\subsection{The case of at most two hyperplanes}\label{sec:twoplanes}

First, let us consider only one halfspace,  $m=1$. As we mention before, this case has been considered in the literature several times. 
The halfspace,
$$H^+_1=\{\xbf \in \reals^n~|~\abf\cdot \xbf + r_1=a_1x_1+a_2x_2+\cdots+a_nx_n+r_1\geq0\}$$
is an auxiliary plane itself. We get $\starv=\varnothing, I=\varnothing$ and $\|\starv\|=|I|=0$, $A_\varnothing^\varnothing=1$ and $A_{i}^1 =a_i$. 
The good clipping condition (A) automatically holds and (B) is equivalent to $\prod_{t=1}^n  a_t \neq 0$.  Applying these terms to Theorem \ref{thm:maincup} we get a proof of Theorem \ref{thm:monthly}.

Second, let us prove  Corollary \ref{thm:twoplanes}.
Consider the following two hyperplanes,
\begin{align*}
H^+_1 &=\{\xbf \in \reals^n~|~\abf\cdot \xbf + r_1=a_1x_1+a_2x_2+\cdots+a_nx_n+r_1\geq0\},\\
H^+_2 &=\{\xbf \in \reals^n~|~\bbf\cdot \xbf + r_2=b_1x_1+b_2x_2+\cdots+b_nx_n+r_2\geq0\}.
\end{align*}
We see that $\starv$ and $I$ become the empty set or a set of only one element. The former case of the empty set  is the same as the  ``one hyperplane'' case. For the case of $I=[1]=\{1\}$, we put $\starv=\{*({\vbf})\}$ then $\bulletv= [n]\setminus \starv$ and get
$$\prod_{t \in \bulletv} A_{\starv \vee t } ^{ I \vee m}
=\prod_{t \in \bulletv}  \left| \begin{array}{cc} a_{*({\vbf})} & b_{*({\vbf})} \\ a_{t} & b_{t} \\ \end{array} \right|.$$
Applying Theorem \ref{thm:mainvee} to these, we obtain Corollary \ref{thm:twoplanes}.
Here, the good clipping conditions are 
\begin{enumerate}[(A)]
	\item $F^0 \cap H_1 \cap H_2^{+}=\varnothing$,
	\item  
	$\prod\limits_{t=1}^n  b_t  \prod\limits_{t \in \bulletv}
	\left| 
	\begin{array}{cc} 
	a_{*({\vbf})} & b_{*({\vbf})} \\ a_{t} & b_{t} \\ 
	\end{array} 
	\right| \neq 0 ~~$ for $\vbf \in F^1 \cap H_1 \cap H_2^{+}$.
\end{enumerate}

\subsection{The case of three hyperplanes.}
Let us consider three halfspaces, $m=3$. 
We formulate this case in a similar fashion to the one or two hyperplane cases, in particular, which is important to derive several identities in Appendix \ref{app:isoceles}.

\begin{corollary}\label{cor:threehyp}
	The volume of the standard unit hypercube
	$[0,1]^n$ intersecting
	the three halfspaces 
	\begin{align*}
	H^+_1 =& \{\xbf \in \reals^n ~|~\abf\cdot \xbf + r_1=a_1x_1+a_2x_2+\cdots+a_nx_n+r_1\geq0\},\\ 
	H^+_2 =& \{\xbf \in \reals^n ~|~ \bbf\cdot \xbf + r_2=b_1x_1+b_2x_2+\cdots+b_nx_n+r_2\geq0\},\\ 
	H^+_3 =& \{\xbf \in \reals^n ~|~ \cbf \cdot \xbf + r_3=c_1x_1+c_2x_2+\cdots+c_nx_n+r_3\geq0\},
	\end{align*}	
	satisfying good clipping conditions, is
	\begin{align*}
	\vol([0,1]^n \cap H_1^{+} \cap H_2^{+} \cap H_3^{+})
	&=
	\sum_{{\vbf}\in F^0 \cap H_1^{+} \cap H_2^{+} \cap H_3^{+}}
	\frac{(-1)^{|\zerov|}g_3({\vbf})^n}{n! \prod_{t\in[n]} c_t}\\
	&-\sum_{{\vbf}\in F^1 \cap H_1 \cap H_2^{+}\cap H_3^{+}}
	\frac{(-1)^{|\zerov|}~sgn(a_{*({\vbf})})~a_{*({\vbf})}^{n-1}~g_3({\vbf})^n}{n!~ c_{*({\vbf})}\prod\limits_{t\in\bulletv} \left| \begin{array}{cc} a_{*({\vbf})} & c_{*({\vbf})} \\ a_{t} & c_{t} \\ \end{array} \right|}\\
	&-\sum_{{\vbf}\in F^1 \cap H_1^+ \cap H_2 \cap H_3^{+}}
	\frac{(-1)^{|\zerov|}~sgn(b_{*({\vbf})})~b_{*({\vbf})}^{n-1}~g_3({\vbf})^n}{n!~ c_{*({\vbf})}\prod\limits_{t\in\bulletv} \left| \begin{array}{cc} b_{*({\vbf})} & c_{*({\vbf})} \\ 	
		b_{t} & c_{t} \\ \end{array} \right|}
	\intertext{\vspace{-0.7em}$$
		-\sum\limits_{\substack{{\vbf}\in F^2 \cap \\H_1  \cap H_2\cap H_3^{+}}}
		\frac{(-1)^{|\zerov|}~sgn(\left| \begin{array}{cc} a_{*_1({\vbf})} & b_{*_1({\vbf})} \\ a_{*_2({\vbf})} & b_{*_2({\vbf})} \\  \end{array} \right|)~\left| \begin{array}{cc} a_{*_1({\vbf})} & b_{*_1({\vbf})} \\ a_{*_2({\vbf})} & b_{*_2({\vbf})}  \\  \end{array} \right|^{n-1}~g_3({\vbf})^n}{n!~ \left| \begin{array}{cc} a_{*_1({\vbf})} & c_{*_1({\vbf})} \\ a_{*_2({\vbf})} & c_{*_2({\vbf})} \\ \end{array} \right| \left| \begin{array}{cc} b_{*_1({\vbf})} & c_{*_1({\vbf})} \\ b_{*_2({\vbf})} & c_{*_2({\vbf})} \\ \end{array} \right|\prod\limits_{t\in\bulletv} \left| \begin{array}{ccc} a_{*_1({\vbf})} & b_{*_1({\vbf})} & c_{*_1({\vbf})} \\ a_{*_2({\vbf})} & b_{*_2({\vbf})} & c_{*_2({\vbf})} \\ a_{t} & b_t & c_{t} \\ \end{array} \right|}.
		$$
		\vspace{-2.5em}
	}
	\end{align*}
\end{corollary}
\begin{proof}
	For each  vertex, $|\starv| = |I| = 0, 1$, or $2$. The former two cases are the same as the case of fewer than two hyperplanes. Let us consider the $I=[2]$ case.
	Recall $\bulletv=[n]\setminus \starv$ and  $\starv=\{*_1({\vbf}),*_2({\vbf})\}$ then
	
	$$
	\prod\limits_{t\in I} A_{\starv}^{I{\vee m\setminus t}}
	=\left| \begin{array}{cc} a_{*_1({\vbf})} & c_{*_1({\vbf})} \\ a_{*_2({\vbf})} & c_{*_2({\vbf})} \\ \end{array} \right| \left| \begin{array}{cc} b_{*_1({\vbf})} & c_{*_1({\vbf})} \\ b_{*_2({\vbf})} & c_{*_2({\vbf})} \\ \end{array} \right|$$
	and
	\begin{align*}
	\prod_{t \in \bulletv } 
	A_{\starv\vee t }^{I\vee m}
	&=\prod_{t\in \bulletv} \left| \begin{array}{ccc} a_{*_1({\vbf})} & b_{*_1({\vbf})} & c_{*_1({\vbf})} \\ a_{*_2({\vbf})} & b_{*_2({\vbf})} & c_{*_2({\vbf})} \\ a_{t} & b_t & c_{t} \\ \end{array} \right|. 
	\end{align*}
	
\end{proof}

\subsection{Examples of calculations }
We show two examples of calculations  using Corollary \ref{thm:twoplanes} and Corollary \ref{cor:threehyp}. In particular  the following examples have several non-simple vertices. But we can apply our formulas to them because  all non-simple vertices lie in the auxiliary hyperplane.
\begin{example}	
	Let us calculate the volume of the clipped hypercube
	$[0,1]^3$ {which intersects}
	the following two halfspaces,
	\begin{align*}
	H_1^{+} &=\{\xbf\in \reals^n ~|~\abf\cdot \xbf+\frac{1}{2}=-x_1+x_2+\frac{1}{2}\geq0\},\\ H_2^{+} &=\{\xbf\in \reals^n ~|~\bbf \cdot \xbf+3=-x_1-2x_2-x_3+3\geq0\},
	\end{align*}

	Let us find the vertices of the clipped hypercube. 
	There are five vertices in $ F^0 \cap H_1^{+} \cap H_2^{+}$: $${\vbf}_1=(0,0,0),~ {\vbf}_2=(0,0,1),~ {\vbf}_3=(0,1,0),~ {\vbf}_4=(0,1,1),~ {\vbf}_5=(1,1,0)$$
	and four vertices   in $F^1 \cap H_1 \cap H_2^{+}$: $${\vbf}_6=(\frac 12,0,0),~ {\vbf}_7=(\frac 12,0,1),~ {\vbf}_8=(1,\frac 12,0),~ {\vbf}_9=(1,\frac 12,1).$$ 
	Among those vertices,  ${\vbf}_4, {\vbf}_5$ and  ${\vbf}_9$ lie on $H_2$ and we don't need to worry about these vertices. We can check that the good clipping conditions hold. 
	We calculate the values $N_{{\vbf}_i}$ for $i=1,2,3,6,7,8$ by Corollary \ref{thm:twoplanes}.
	For example, we have
	$$N_{{\vbf}_6}=-\frac{(-1)^{2}~sgn(-1)~(-1)^2~g_2(\frac 12,0,0)^3}{3!~ (-1) \left|
		\begin{array}{cc}
		-1 & -1 \\
		1 & -2 \\
		\end{array}
		\right| \left|
		\begin{array}{cc}
		-1 & -1 \\
		0 & -1 \\
		\end{array}
		\right| }=-\frac{(\frac{5}{2})^3}{6\times3\times1}=-\frac{125}{144}.$$
	Therefore we get
	$$\begin{aligned}
	\vol([0,1]^3 \cap H_1^{+} \cap H_2^{+})
	=&\sum_{{\vbf}\in\{{\vbf}_1,{\vbf}_2,{\vbf}_3\}}
	\frac{(-1)^{|\zerov|}g_2({\vbf})^3}{3! \prod_{t=1}^3 b_t}\\
	-&\sum_{{\vbf}\in\{{\vbf}_6,{\vbf}_7,{\vbf}_8\}}
	\frac{(-1)^{|{\zerov}|}~a_{*({\vbf})}^{3}~g_2({\vbf})^3}{3!~|a_{*({\vbf})}|~ b_{*({\vbf})}\prod\limits_{t\in [3]\setminus *({\vbf})} \left| \begin{array}{cc} a_{*({\vbf})} & b_{*({\vbf})} \\ a_{t} & b_{t} \\ \end{array} \right|}\\
	=&\frac{9}{4}-\frac{2}{3}-\frac{1}{12}-\frac{125}{144}+\frac{27}{144}-\frac{1}{36}\\=&\frac{19}{24}.\end{aligned}
	$$
	
	\begin{remark}
		Note that the polyhedron $[0,1]^3 \cap H_1^{+} \cap H_2^{+}$ has three non-simple vertices ${\vbf}_4, {\vbf}_5, {\vbf}_9$ but all these are degenerate into $H_2$ so we can apply the formula of Corollary \ref{thm:twoplanes}. Notice that if one changes the roles of the two halfspaces then one cannot apply the formula,  because there are non-degenerate non-simple vertices and it violates the good clipping conditions.	
	\end{remark}
\end{example}	

\begin{example}
	Let us calculate the volume of the region of
	$[0,1]^3$ {that intersects}
	the three halfspaces,
	\begin{align*}
	H_1^{+}&=\{x|~-x_1+x_2+\tfrac{1}{2}\geq0\},\\ 
	H_2^{+}&=\{x|~x_3- \tfrac{1}{2}\geq0\}, \\ 
	H_3^{+}&=\{x|~-x_1-2x_2-x_3+3\geq0\}.	
	\end{align*}
	
	The three halfspaces satisfy the good clipping conditions and  we can apply Corollary~\ref{cor:threehyp}.
	Let us find vertices for each $I \subset [3-1]$, i.e. for $I = \varnothing, \{1\}, \{2\}$ and $\{1,2\}$:	
	\begin{align*}
	F^0 \cap H_1^{+} \cap H_2^{+} \cap H_3^{+} ~~&:~~ {\vbf}_1=(0,0,1),~ {\vbf}_2=(0,1,1),\\ 
	F^1 \cap H_1~ \cap H_2^{+} \cap H_3^{+} ~~&:~~  {\vbf}_3=(\tfrac 12,0,1),~ {\vbf}_4=(1,\tfrac 12,1), \\
	F^1 \cap H_1^{+} \cap H_2~ \cap H_3^{+} ~~&:~~ {\vbf}_5=(0,0,\tfrac 12),~ {\vbf}_6=(0,1,\tfrac 12), \\ 
	F^2 \cap H_1~ \cap H_2~ \cap H_3^{+} ~~&:~~ {\vbf}_7=(\tfrac 12,0,\tfrac 12),~ {\vbf}_8=(1,\tfrac 12,\tfrac 12).
	\end{align*}
	Let us check that ${\vbf}_2$ and ${\vbf}_4$ lie on $H_3$ and these vertices are degenerate vertices. Note that  there are two more degenerate simple vertices $\vbf_9=(1,\tfrac{3}{4},\tfrac{1}{2})$ and $\vbf_{10}=(\tfrac{1}{2},1,\tfrac{1}{2})$  which are excluded from the summation.  
	In summary, the polyhedron $[0,1]^3 \cap H_1^{+} \cap H_2^{+}\cap H_3^{+}$ has ten vertices with two non-simple vertices ${\vbf}_2$ and  ${\vbf}_4$ among them. After applying  Corollary {\ref{cor:threehyp}} to these, 
	we obtain 
	$$N_{\vbf_1}=-\frac{2}{3},~ N_{{\vbf}_3}=\frac{3}{16},~ N_{{\vbf}_5}=\frac{125}{96},~ N_{{\vbf}_6}=-\frac{1}{96},~ N_{{\vbf}_7}=-\frac{4}{9},~\text{ and }~ N_{{\vbf}_8}=-\frac{1}{288}.$$
	
	Therefore we obtain the volume as follows,
	$$\begin{aligned}
	\vol([0,1]^3 \cap H_1^{+} \cap H_2^{+} \cap H_3^{+})
	=&\sum_{{\vbf}=\vbf_1}N_{\vbf} +\sum_{{\vbf}={\vbf}_3}N_{\vbf} +\sum_{{\vbf}\in\{{\vbf}_5,{\vbf}_6\}}N_{\vbf} +\sum_{{\vbf}\in\{{\vbf}_7,{\vbf}_8\}}N_{\vbf}\\
	=&-\frac{2}{3}+\frac{3}{16}+\frac{125}{96}-\frac{1}{96}-\frac{4}{9}-\frac{1}{288}\\=&\frac{35}{96}.\end{aligned}
	$$
	
\end{example}

\section{Combinatorial identities from clipping hypercubes}\label{sec:combidentity}
\subsection{From polytopes to identities}
Let us describe a general method to produce a combinatorial identity  from a polytope volume. This is a simple observation that the resulting volume is independent of the choice of an auxiliary plane. 
Recall the volume expression of Theorem \ref{thm:mainLawrence} and theorems in Section \ref{sec:severalversion} and let us assume that  we already know the volume of a clipped hypercube $P = [0,1]^n \cap {H_1^+} \cap \dots \cap {H_{m-1}^+}$.  Let us cut  $P$ into two pieces {one more time} as
$$P_+ = P \cap H_m^+ ~~~\text{ and }~~ P_- = P \cap H_m^-$$ 
by the auxiliary hyperplane $H_{m}=\{a_1 x_1 + \dots + a_n x_n +y = 0\}.$
No matter how we take $H_m$, the union of two pieces should be $P$ and  $$\vol(P_+) + \vol(P_-)=\vol(P). $$  
The known volume is constant and is expressed in terms of the  free variables $a_1,a_2,\dots,a_n$ and $y$ of the coefficients of $H_m$. They will produce an algebraic identity.
\begin{remark}
	Only hyperplanes satisfying  good clipping conditions actually produce an algebraic identity since the volume formula makes sense only for this case. However, as we discussed in Section \ref{sec:epsilonpert}, we can use a limiting argument since the volume function is continuous function and the good clipping conditions are  open conditions. Therefore, whether a good clipping condition is satisfied or not, the resulting algebraic identity holds  as long as the expression is valid. 
\end{remark}
\begin{remark}
	The volume formulas are homogeneous for $a_1,a_2,\dots a_n$ because they are composed of  homogeneous polynomials which are determinants of  matrices with  one column vector of  indeterminate  $\abf = (a_1,a_2,\dots, a_n )$.  
\end{remark} 
Let  us see the most simple case which is a direct consequence of Theorem \ref{thm:monthly}.

\begin{corollary}\label{cor:1hyperplane} Let $\abf = (a_1,\dots,a_n) \in \reals^n$ and $y \in \reals$. 
	Then 
	$$	\sum_{ {\vbf} \in F^0}
	{(-1)^{|\zerov|} (\abf \cdot {\vbf}+y)^n}= {n! a_1a_2\cdots a_n}.
	$$
\end{corollary}
\begin{proof} 
	
	For a hyperplane $$H_1=\{\xbf ~|~ g(\xbf):= \abf \cdot \xbf +y = a_1x_1+ \cdots + a_n x_n  + y  =0\},$$ 
	{we get}
	$$\vol ( [0,1]^n \cap {H_1^-})	= \sum_{ {\vbf} \in F^0 \cap {H_1^-}}
	\frac{(-1)^{|\zerov|} g({\vbf})^n} {n! \prod_{t=1}^n  a_t}$$ and $$\vol ( [0,1]^n)=\vol ( [0,1]^n \cap {H_1^+})+\vol ( [0,1]^n \cap {H_1^-}),$$
	so we have $$1=\sum_{ {\vbf} \in F^0}
	\frac{(-1)^{|\zerov|} g({\vbf})^n} {n! \prod_{t=1}^n  a_t}.$$  
	Finally, we get 
	\begin{align*}
	\sum_{ {\vbf} \in F^0}
	{(-1)^{|\zerov|} (\abf \cdot {\vbf}+y)^n} &= {n! \prod_{t=1}^n  a_t}.
	\end{align*}
\end{proof}
Let us consider the summation over  $\vbf \in F^0 =   \{0,1\}^n$. 
We  can replace this summation by  {
$1 \leq {t_1},\dots {t_i}  \leq n$ } as regarding {$\onev =  \{{t_1},\dots,{t_i}\}$ }and hence prove Theorem \ref{thm:generalizedRuiz}.   Note that $a_i$ should be non-zero when applying  the volume formula but the resulting identity has no such constraint by continuity, as we remarked above.

Essentially, whenever we take a  polytope, we can find a corresponding combinatorial identity if we have a concrete volume formula.
Hence we can expect this kind of 
$$\{~polytopes~\} \longrightarrow \{~combinatorial~~ identities~\} $$ correspondence has a structural property. 
At this stage, it seems  to be somewhat vague  to investigate the resulting identities from general convex polytopes. 
Here, we present a few cases. We give the case of a clipped hypercube by a symmetric hyperplane with full generality in the next section, which produces the  interesting identity in Theorem \ref{thm:finalcom}.  
We treat several more examples of resulting identities in the Appendix.

\subsection{Symmetric arrangements of  hyperplanes}
When we see  the identity of Theorem \ref{thm:generalizedRuiz}, we can observe that this is a symmetric function of the $n$-variables $a_1, \dots a_n$. 
This property comes from the fact that the  polytope itself under consideration is symmetric, i.e. we took a symmetric arrangement of  hyperplanes, where the term  \emph{symmetric} means that hyperplanes except the auxiliary hyperplane are invariant under exchange of coordinate axes of $\reals^n$. 

Probably the second easiest example of a symmetric arrangement is 
$${H_1^+}=\{\xbf \in \reals^n ~|~ -x_1-x_2-\dots-x_n + 1 \geq 0 \}.$$
We use the auxiliary {hyperplane} 
$$H_2=\{\xbf \in \reals^n ~|~ a_1 x_1+\dots a_n x_n + y {=} 0\}$$
in the formula of the $m=2$ case  and obtain the  identity of Theorem \ref{thm:clippedsimplex}. We remark that the identity of Theorem \ref{thm:clippedsimplex}  is  a direct consequence of  Proposition 1 in \cite{cho_volume_2001}.

Note that $H_1$ and $H_2$ violate the good clipping conditions. So we use  the $\epsilon$-perturbation proposed in Section \ref{sec:goodclipping} when applying the volume formula. 
\begin{proof}[Proof of Theorem \ref{thm:clippedsimplex}] For sufficiently small $\epsilon>0$, the polytope $[0,1]^n \cap H_1^+$ with  
	$$H_1^{+}=\{\xbf|~-x_1-x_2-\cdots -x_n+1-\epsilon\geq0\}$$
	 satisfies the good clipping conditions. We have
	$$\begin{aligned}
	\vol ( [0,1]^n \cap H_1^+)&=\vol ( [0,1]^n \cap H_1^+\cap H_2^+)+\vol ( [0,1]^n \cap H_1^+\cap H_2^-)\\
	&=\sum_{{\vbf}\in F^0 \cap H_1^{+} \cap H_2^{+}} - \sum_{{\vbf}\in F^1 \cap H_1 \cap H_2^{+}}+\sum_{{\vbf}\in F^0 \cap H_1^{+} \cap H_2^{-}} - \sum_{{\vbf}\in F^1 \cap H_1 \cap H_2^{-}}\\
	&=\sum_{{\vbf}\in F^0 \cap H_1^{+}} \frac{(-1)^{|\zerov|}g_2({\vbf})^n}{n! \prod_{t=1}^n a_t}\\
	-&\sum_{{\vbf}\in F^1 \cap H_1 }
	\frac{(-1)^{|\zerov|}~(-1)^{n}~g_2({\vbf})^n}{n!~ a_{*({\vbf})}\prod_{t=1,t\neq *({\vbf})}^n \left| \begin{array}{cc} -1 & a_{*({\vbf})} \\ -1 & a_{t} \\ \end{array} \right|}.
	\end{aligned}$$
	
	Then there are $n+1$ vertices of 
	\begin{alignat*}{2}
	{\vbf} &=(0,0,\ldots,0) &&~~~\in~ F^0 \cap H_1^{+}, \\
	\vbf_i &= (1-\epsilon)\ebf_i &&~~~\in~  F^1 \cap H_1 ~~\text{ for }~  i \in [n]. 
	\end{alignat*}
	
	Hence we obtain
	$$\frac{(1-\epsilon)^n}{n!}=\frac{(-1)^n y^n}{n!a_1a_2\cdots a_n}-\sum_{i=1}^n \frac{(-1)^{n-1}(-1)^n(a_i(1-\epsilon)+y)^n}{n!a_i \prod_{t=1,t\neq i}^n (a_i - a_t)}.$$
	The result now follows by taking $\epsilon\rightarrow 0$ and simplifying. 	
\end{proof}
We next consider the following {half space},
$${H_1^+}=\{\xbf \in \reals^n ~|~ -x_1-x_2-\dots-x_n + 2 \geq 0 \}.$$
By the same  $\epsilon$-perturbation taking $2-\epsilon$ instead of $2$, we get the following identity,
\begin{align*}
&\frac{ y^n}{a_1a_2\cdots a_n} ~~-~~ \sum_{i=1}^n \frac{(y+a_i)^n}{a_1a_2\cdots a_n}+\\
&\sum_{1\leq t_1<t_2\leq n}\sum_{i=1}^2 \frac {(y+a_{t_1}+a_{t_2})^n}{a_{t_i}(a_1-a_{t_i})(a_2-a_{t_i})\cdots (a_n-a_{t_i})}
= (-1)^n (2^n-n).
\end{align*}

We  consider all possible symmetric arrangements of only one hyperplane.  Then  all of the linear coefficients of $H_1$ should be the same. So 
it is reasonable to think about
$${H_1^+}=\{\xbf \in \reals^n ~|~ -x_1-x_2-\dots-x_n + l \geq 0 \}
~~~\text{ for }~~ l \in [n].$$
This  gives the following theorem which is nothing but a different form  of Theorem \ref{thm:finalcom}.  
\begin{theorem}\label{thm:finalid}  For an integer $l \in [n]$  and non-zero distinct real numbers $a_1, a_2,\ldots, a_n$ and $y \in \reals$, 
	\begin{equation*}\label{e4} 		
	\begin{aligned}
	\frac{y^n}{a_1a_2\cdots a_n} +&~ \sum_{i=1}^{l-1} \sum_{1\leq t_1<t_2<\cdots<t_i\leq n} \frac {(-1)^{i}(y+a_{t_1}+a_{t_2}+\cdots +a_{t_i})^n}{a_1a_2\cdots a_n}\\
	+&~(-1)^{l}\sum_{1\leq t_1<\cdots <t_l\leq n}\sum_{i=1}^l \frac {(y+a_{t_1}+\cdots+a_{t_l})^n}{a_{t_i}(a_1-a_{t_i})(a_2-a_{t_i})\cdots (a_n-a_{t_i})} \\
	=&~~ \sum_{i=0}^{l-1} (-1)^{n-i} \binom{n}{i}(l-i)^n
	\end{aligned}\end{equation*}
	
	or equivalently
	$$\begin{aligned}&\sum_{i=1}^{l-1} \sum_{1\leq t_1<t_2<\cdots<t_i\leq n}  {(-1)^{i}(a_{t_1}+a_{t_2}+\cdots +a_{t_i})^k}\\
	&+(-1)^{l}\sum_{1\leq t_1<\cdots <t_l\leq n}\sum_{i=1}^l (\prod_{j=1,j\ne t_i}^n\frac {a_j}{a_j-a_{t_i}})    {(a_{t_1}+\cdots+a_{t_l})^k}\\
	&=                   \begin{cases}
	-1, & \textup{if }\;\; k=0 \\
	0 ,  & \textup{if }\;\; k=1,2,\ldots,n-1 \\
	\prod_{i=1}^n{a_i}\sum_{i=0}^{l-1} (-1)^{n-i} \binom{n}{i}(l-i)^n ,  & \textup{if }\;\; k=n
	\end{cases}
	\end{aligned}$$
	or,  using set-notation with $A=\{a_1, a_2,\ldots, a_n\}$, 
	\begin{equation*}\label{e5}		
	\begin{aligned} &\sum_{\substack{I\subset A\\|I|<l}}  {(-1)^{|I|}(y+\|I\|)^n}
	+\sum_{\substack{I\subset A\\|I|=l}} (-1)^{l}(y+\|I\|)^n (\sum_{a\in I} \prod_{b\in A\setminus a}\frac {b}{b-a})\\ 
	&=A! \sum_{i=0}^{l-1} (-1)^{n-i} \binom{n}{i}(l-i)^n.
	\end{aligned}
	\end{equation*}
\end{theorem}
One needs to be careful here; there is a difference between Theorem \ref{thm:finalcom} and   the second form above in the case $k=0$ because $\|\varnothing \|^0 = 0^0= 1$.

\begin{proof} We use  $\epsilon$-perturbation  replacing $l$ by $l-\epsilon$ for $H_1$, then
	$$\begin{aligned}
	\vol ( [0,1]^n \cap H_1^+)&=\sum_{{\vbf}\in F^0 \cap H_1^{+}} \frac{(-1)^{|\zerov|}(\abf\cdot {\vbf}+y)^n}{n!~a_1a_2\cdots a_n}\\
	-&\sum_{{\vbf}\in F^1 \cap H_1 }
	\frac{(-1)^{n-l}~(-1)^n~(\abf\cdot {\vbf}+y)^n}{n!~a_{*({\vbf})}\prod_{t=1,t\neq *({\vbf})} (a_{*({\vbf})}-a_t)}.
	\end{aligned}$$
	
	For $ F^0 \cap H_1^{+}$,	there are a total of $\binom{n}{0}+\binom{n}{1}+\cdots+\binom{n}{l-1}$  vertices, i.e. there are $l$  families of vertices with respect to the sum 
	of {the} coordinate values of ${\vbf}$. 
	
	Also, for $F^1 \cap H_1$,
	there are $l\binom{n}{l}$ vertices whose coordinate values are $l-1$   ones, a unique $1-\epsilon$, and $n-l$  zeros.
	Hence we obtain 
	\begin{align*}
	&\sum_{{\vbf}\in F^0 \cap H_1^{+}} \frac{(-1)^{|\zerov|}(\abf\cdot {\vbf}+y)^n}{n!~a_1a_2\cdots a_n}\\
	&~~~~~~~~~~~=\frac{(-1)^n y^n}{n!~a_1a_2\cdots a_n} + \sum_{i=1}^{l-1} \sum_{1\leq t_1<\cdots<t_i\leq n} \frac {(-1)^{n-i}(y+a_{t_1}+\cdots +a_{t_i})^n}{n!~a_1a_2\cdots a_n}\
	\intertext{and}
	&\sum_{{\vbf}\in F^1 \cap H_1 }
	\frac{(-1)^{l+1}~(\abf\cdot {\vbf}+y)^n}{n!~a_{*({\vbf})}\prod_{t=1,t\neq *({\vbf})} (a_{*({\vbf})}-a_t)}\\
	&~~~~~~~~~~~=(-1)^{l+1}\sum_{1\leq t_1<\cdots <t_l\leq n}\sum_{i=1}^l \frac {(y+a_{t_1}+\cdots+a_{t_i}(1-\epsilon) +\cdots +a_{t_l})^n}{n!~a_{t_i}(a_{t_i}-a_1)(a_{t_i}-a_2)\cdots (a_{t_i}-a_n)}.
	\end{align*}
	We compute the volume of the clipped hypercube using Theorem \ref{thm:monthly},
	$$\begin{aligned}
	\vol ( [0,1]^n \cap H_1^+)=&\frac{(-1)^n(l-\epsilon)^n}{n!(-1)^n}+\frac{(-1)^{n-1} \binom{n}{1}(l-1-\epsilon)^n}{n!(-1)^n}+\cdots\\&+\frac{(-1)^{n-(l-1)} \binom{n}{l-1}(l-(l-1)-\epsilon)^n}{n!(-1)^n}\\
	=&\sum_{i=0}^{l-1} (-1)^{i} \binom{n}{i}\frac{(l-i-\epsilon)^n}{n!}.
	\end{aligned}$$
	By taking $\epsilon\rightarrow 0$, we obtain the result. 	
\end{proof}
\begin{remark}
	If we take $l$ to be a non-integer real number, we get a slightly different identity,  obtained  by rescaling  variables  from the result of Theorem \ref{thm:finalid}. If the vertex configuration is preserved under changing hyperplanes, the resulting identity is essentially the same as the previous one.   	
\end{remark}

%
%

\begin{appendix}

	\renewcommand{\thesection}{\hspace{-0.3em}}
	
	\section{Several clipped hypercube identities}\label{sec:appendix}

	For simplicity,  we do not use $m$ for the number of hyperplanes in the appendix section and  use $\obf_n$ instead of $(0,0,\ldots,0)$ in $\reals^n$. 
	
	\renewcommand{\thesubsection}{\Alph{subsection}}
	
	\subsection{Symmetric truncated hypercube}
	Let us consider $n+1$ hyperplanes
	\begin{align*}
	H_1 &=\{\xbf ~|~-x_1+x_2+\cdots+ x_n+1-d=0\}, \\
	H_2 &=\{ \xbf ~|~ x_1-x_2+\cdots+x_n+1-d=0\},\\
	& \hspace{7em}\vdots\\
	H_n &=\{\xbf ~|~x_1+x_2+\cdots- x_n+1-d=0\}, \\
	H_{n+1} &=\{\xbf ~|~a_1x_1+a_2x_2+\cdots+a_n x_n+y=0\}. \\
	\end{align*}
	Then the volume is the following,
	\begin{align*}
	&\vol([0,1]^n \cap H_1^{+} \cap \cdots  \cap H_{n}^{+})\\
	&~~~~~~~=
	\vol([0,1]^n \cap H_1^{+} \cap \cdots  \cap H_{n}^{+}\cap H_{n+1}^{+})+ \vol([0,1]^n \cap H_1^{+} \cap \cdots  \cap H_{n}^{+}\cap H_{n+1}^{-})\\
	&~~~~~~~= 1-n\times \frac{d^n}{n!}.
	\end{align*}
	
	Note that these hyperplanes do not intersect  each other in $[0,1]^n$ under the condition $0<d<1$.
	Corollary \ref{thm:twoplanes} essentially suffices to compute the volume. 
	
	We can check that there are  three kinds of vertices,
	\begin{align*}
	|I|=0 :&~~~F^0\setminus \{\ebf_1,\ebf_2,\dots,\ebf_n\},\\
	|I|=1 :&~ (1-d)\ebf_i ~\text{ for }~ i = 1,2,\dots,n ,\\
	&~ \ebf_i+d\ebf_j ~~\text{ for }~ 1\le i\ne j \le n.
	\end{align*}
	The resulting identity  is
	\begin{align*}
	&\sum_{i=1}^{n}\frac{ (y+a_i)^n}{a_1a_2\cdots a_n} - \sum_{i=1}^{n} \frac {(y+a_{i}(1-d))^n}{a_i\prod_{j=1,j\ne i}^n (a_j+a_i)}\\
	&-\sum_{1\leq i\ne j \leq n} \frac {(y+a_{i}+a_{j}d)^n}{a_{j}(a_{i}+a_j)\prod_{t=1,t\ne i,j}^n(a_t-a_j)} =(-1)^{n+1} n d^n.
	\end{align*}

	\subsection{Hyperprism : $n\text{-simplex }\times [0,1]^m$}
	Let us consider the following two hyperplanes
	\begin{align*}
	H_1 &=\{\xbf ~|~-x_1-x_2-\cdots- x_n+1-\epsilon=0\}, \\
	H_2 &=\{\xbf ~|~a_1x_1+\cdots+a_n x_n+b_1x_{n+1}+\cdots+b_m x_{n+m}+y=0\}. \\
	\end{align*}
	The resulting volume taking $\epsilon\rightarrow0$ is the following.
	\begin{alignat*}{2}
	\vol ([0,1]^{n+m} \cap H_1^{+})~&&=&~\vol([0,1]^{n+m} \cap H_1^{+} \cap  H_{2}^{+})+ \vol([0,1]^{n+m} \cap H_1^{+} \cap  H_{2}^{-})\\
	&&=&~~\frac{1}{n!}.
	\end{alignat*}
	
	We can check that there are several kinds of vertices 
	\begin{align*}
	|I|=0: &&  &\obf_{n+m},\\
	&& &\ebf_{n+1},\ebf_{n+2},\ldots,\ebf_{n+m},\\
	&& &\ebf_{n+1}+\ebf_{n+2},\ebf_{n+1}+\ebf_{n+3},\ldots,\ebf_{n+m-1}+ \ebf_{n+m},\\
	&& &\hspace{7em}\vdots\\
	&& & \ebf_{n+1}+\ebf_{n+2}+\cdots+\ebf_{n+m},\\
	\intertext{and for  $i \in [n]$, }
	|I|=1: &&  &(1-\epsilon)\ebf_i+\obf_{n+m},\\
	&& &(1-\epsilon)\ebf_i+ \{\ebf_{n+1},\ldots,~ \ebf_{n+m}\},\\
	&& &(1-\epsilon)\ebf_i+ \{ \ebf_{n+1}+\ebf_{n+2},\ebf_{n+1}+\ebf_{n+3},\ldots,~\ebf_{n+m-1}+ \ebf_{n+m} \},\\
	&& &\hspace{7em}\vdots\\
	&& & (1-\epsilon)\ebf_i+\ebf_{n+1}+\ebf_{n+2}+\cdots+\ebf_{n+m}.
	\end{align*}

	The resulting identity is
	\begin{align*}
	&\sum_{I\subset\{b_1,\ldots,b_m\}} {(-1)^{|I|} (y+\|I\|)^{n+m}} \\
	&~~~~~~~~~~~+ \sum_{i=1}^{n} \sum_{I\subset\{b_1,\ldots,b_m\}} {(-1)^{|I|+1}(\prod_{j=1,j\ne i}^n \frac{a_j}{a_j-a_i})(y+a_{i}+\|I\|)^{n+m}}\\
	&~~~~~~~~~~~=(-1)^{n+m}\frac {(n+m)!}{n!} {a_1\cdots a_n b_1\cdots b_m}.
	\end{align*}

	\subsection{Isosceles $n$-simplex} \label{app:isoceles}
	Let us consider the following three hyperplanes
	\begin{align*}
	H_1 &=\{\xbf ~|~-x_1-x_2-\cdots- x_n+1-\epsilon=0\}, \\
	H_2 &=\{ \xbf ~|~ x_1-x_2-\cdots-x_n-\epsilon=0\},\\
	H_3 &=\{\xbf ~|~a_1x_1+a_2x_2+\cdots+a_n x_n+y=0\}.
	\end{align*}
	The resulting volume  taking $\epsilon\rightarrow0$ is the following.
	\begin{align*}
	&\vol([0,1]^n \cap H_1^{+} \cap H_2^{+})\\
	&=\vol([0,1]^n \cap H_1^{+} \cap H_2^{+}\cap H_3^{+})+ \vol([0,1]^n \cap H_1^{+} \cap  H_2^{+}\cap H_3^{-})\\
	&=\frac{1}{n!2^{n-1}}.
	\end{align*}
	We can check that there are two kinds of vertices 
	\begin{alignat*}{2}
	|I|=1 :&~~~\epsilon\ebf_1,(1-\epsilon)\ebf_1 ,\\
	|I|=2 :&~~~\frac{1}{2}\ebf_1+(\frac{1}{2}-\epsilon)\ebf_i ~\text{ for }~  2\le i \le n.
	\end{alignat*}
	This case needs Corollary \ref{cor:threehyp} for the three hyperplane case.
	
	The resulting identity is
	\begin{align*}
	&\frac{ y^n}{a_1(a_2+a_1)(a_3+a_1)\cdots(a_n+a_1)}-\frac{(y+a_1)^n}{a_1(a_2-a_1)(a_3-a_1)\cdots(a_n-a_1)}\\
	&-2\sum_{i=2}^{n} \frac {(y+\frac{a_1}{2}+\frac{a_i}{2})^n}
	{(a_1+a_i)(a_1-a_i)(a_2-a_i)\cdots(a_n-a_i)} =(-1)^{n} 2^{1-n}.
	\end{align*}

	\subsection{Trapezoidal polytope}
	Let us consider the following two hyperplanes
	\begin{align*}
	H_1 &=\{\xbf ~|~-\frac{x_1}{2}-\frac{x_2}{2}-\cdots--\frac{x_n}{2}-x_{n+1}-x_{n+2}-\cdots- x_{n+m}+1-\epsilon=0\}, \\
	H_2 &=\{\xbf ~|~a_1x_1+\cdots+a_n x_n+b_1x_{n+1}+\cdots+b_m x_{n+m}+y=0\}.
	\end{align*}
	The resulting volume  taking $\epsilon\rightarrow0$ is the following.
	\begin{alignat*}{2}
	\vol([0,1]^{n+m} \cap H_1^{+})&&=&~~\vol([0,1]^{n+m} \cap H_1^{+} \cap  H_{2}^{+})+ \vol([0,1]^{n+m} \cap H_1^{+} \cap  H_{2}^{-})\\
	&&=&~~\frac{2^n-n2^{-m}}{(n+m)!}.
	\end{alignat*}
	
	We can check that there are four kinds of vertices:
	\begin{align*}
	|I|=0 : ~~~& \obf_{n+m},~ \ebf_1,~ \ebf_2,\ldots,~\ebf_n, \\
	|I|=1: ~~~& (1-\epsilon)\ebf_{n+i} \;~~~~~~~~\text{ for }~ 1\le i \le m ,\\
	&\ebf_i+(\frac{1}{2}-\epsilon)\ebf_{n+j}  ~~\text{ for }~ 1\le i \le n, 1\le j \le m,\\
	&\ebf_i+(1-2\epsilon)\ebf_j  \,~~~~\text{ for }~ 1\le i\ne j \le n.
	\end{align*}
	The resulting identity is
	\begin{align*}
	&\frac{(-1)^{n+m} y^{n+m}}{a_1\cdots a_n b_1\cdots b_m} +\sum_{i=1}^n \frac{(-1)^{n+m-1} (y+a_i)^{n+m}}{a_1\cdots a_n b_1\cdots b_m}\\
	&+\sum_{j=1}^{m} \frac {(y+b_j)^{n+m}}{b_j\prod_{s=1}^n (\frac{b_j}2-a_s)\prod_{t=1,t\ne j}^m(b_j-b_t)}\\
	&-\sum_{i=1}^n \sum_{j=1}^{m} \frac {(y+a_i+\frac{b_j}{2})^{n+m}}{b_j\prod_{s=1}^n (\frac{b_j}2-a_s)\prod_{t=1,t\ne j}^m(b_j-b_t)}\\
	&-\sum_{{1\leq i\ne j \leq n}} \frac {(y+a_i+a_j)^{n+m}}{2^m a_i\prod_{s=1,s\ne i}^n (a_i-a_s)\prod_{t=1}^m(a_i-\frac{b_t}2)}\\
	&=2^n-n2^{-m}.
	\end{align*}

\end{appendix}

\bigskip
{\bf Acknowledgment.} 
The authors would like to thank professor Eungchun Cho at Kentucky State University for introducing this problem and also thank   Gabriel C. Drummond-Cole for his helpful comments.  
The authors appreciate several anonymous referees for many valuable  suggestions and corrections.


\end{document}